\documentclass[a4paper]{article}
\usepackage{amsthm}
\usepackage{amsfonts}
\usepackage{amssymb}
\usepackage{amsmath}
\usepackage{graphicx} 
\usepackage{hyperref}
\usepackage{xcolor}
\hypersetup{
	colorlinks,
	linkcolor={blue!50!black},
	citecolor={blue!50!black},
	urlcolor={blue!50!black}
}
\usepackage{cleveref}
\newtheorem{theorem}{Theorem}[section]
\newtheorem{lemma}[theorem]{Lemma}

\newtheorem{proposition}[theorem]{Proposition}
\newtheorem{corollary}[theorem]{Corollary}
\newtheorem{definition}[theorem]{Definition}

\newtheorem{remark}[theorem]{Remark}

\newtheorem{example}[theorem]{Example}
\newtheorem{question}[theorem]{Question}
\newtheorem{problem}[theorem]{Problem}
\newcommand{\N}{\mathbb{N}}

\newcommand{\Z}{\mathbb{Z}}

\newcommand{\M}{\mathfrak{M}}
\newcommand{\G}{\Gamma}
\newcommand{\cay}{\operatorname{Cay}}
\newcommand{\J}{\mathcal{J}}
\renewcommand{\L}{\Lambda}
\newcommand{\e}{\epsilon}
\newcommand{\F}{\mathcal{F}}
\begin{document}
	\title{Translation-like actions by $\mathbb{Z}$, the subgroup membership problem, and Medvedev degrees of effective subshifts}
	\author{Nicanor Carrasco-Vargas}
	\date{}
	\maketitle
	\begin{abstract} 
		We show that every infinite, locally finite, and connected graph admits
		a translation-like action by $\mathbb{Z}$, and that this action can be taken
		to be transitive exactly when the graph has either one or two ends.
		The actions constructed satisfy $d(v,v\ast 1)\leq3$ for every vertex
		$v$. This strengthens a theorem by Brandon Seward.  We also study the effective computability of translation-like actions
		on groups and graphs. We prove that every finitely generated infinite
		group with decidable word problem admits a translation-like action
		by $\Z$ which is computable, and satisfies an extra condition which
		we call decidable orbit membership problem.  As a nontrivial application of our results, we prove that for every
		finitely generated infinite group with decidable word problem, effective
		subshifts attain all $\Pi_{1}^{0}$ Medvedev degrees. This extends
		a classification proved by Joseph Miller for $\Z^{d},$ $d\geq1$. \end{abstract}
	
	\section{Introduction}
		\subsection{Translation-like actions by $\Z$ on locally finite graphs}
		
		A right action $\ast$ of a group $H$ on a metric space $(X,d)$
		is called a \textit{translation-like} \textit{action} if it is \textit{free} (that is, $x\ast h=x$ implies $h=1_{H}$, for $x\in X$, $h\in H$), 
		and for each $h\in H$, the set $\{d(x,x\ast h)|\ x\in X\}\subset\mathbb{R}$
		is bounded. If $G$ is a finitely generated group endowed with the
		left-invariant word metric associated to some finite set of generators,
		then the action of any subgroup $H$ on $G$ by right translations
		$(g,h)\mapsto gh$ is a translation-like action. On the other hand,
		observe that despite the action $H\curvearrowright G$ by left multiplication
		is usually referred to as an action by translations, in general it
		is not translation-like for a left-invariant word metric. 
		
		Following this idea, Kevin Whyte proposed in \cite{whyte_amenability_1999}
		to consider translation-like actions as a generalization of subgroup
		containment, and to replace subgroups by translation-like actions
		in different questions or conjectures about groups and subgroups.
		This was called a geometric reformulation. For example, the von Neumann
		Conjecture asserted that a group is nonamenable if and only if it
		contains a nonabelian free subgroup. Its geometric reformulation asserts
		then that a group is nonamenable if and only if it admits a translation-like
		action by a nonabelian free group. While the conjecture was proven
		to be false \cite{olshanskij_question_1980}, Kevin Whyte proved that
		its geometric reformulation holds \cite{whyte_amenability_1999}. 
		
		One problem left open in \cite{whyte_amenability_1999} was the geometric
		reformulation of Burnside's problem. This problem asked if every finitely
		generated infinite group contains $\Z$ as a subgroup, and was answered
		negatively in \cite{golod_class_1964}. Brandon Seward proved that
		the geometric reformulation of this problem also holds. 
		\begin{theorem}[Geometric Burnside's problem, \cite{seward_burnside_2014}]
			\label{thm:seward} Every finitely generated infinite group admits
			a translation-like action by $\Z$. 
		\end{theorem}
		
		A finitely generated infinite group with two or more ends has a subgroup
		isomorphic to $\Z$, by Stalling's structure theorem. Thus, it is
		the one ended case that makes necessary the use of translation-like
		actions. In order to prove \Cref{thm:seward}, Brandon Seward
		proved a more general graph theoretic result:
		\begin{theorem}[{\cite[Theorem 1.6]{seward_burnside_2014}}]
			\label{thm:sewardgrafos} Let $\G$ be a connected and infinite graph
			whose vertices have uniformly bounded degree. Then $\G$ admits a
			transitive translation-like action by $\Z$ if and only if it is connected
			and has either one or two ends. 
		\end{theorem}
		
This result proves \Cref{thm:seward} for groups with one or
		two ends, and indeed it says more, as the translation-like action
		obtained is transitive. The proof of this result relies strongly on
		the hypothesis of having uniformly bounded degree. Indeed, the uniform
		bound on $d_{\G}(v,v\ast1)$ depends linearly on a uniform bound for
		the degree of the vertices of the graph. Here we strengthen Seward's
		result by weakening the hypothesis to the locally finite case, and
		improving the bound on $d_{\G}(v,v\ast1)$ to 3. 
		
		\begin{theorem}
			\label{thm:t1-translation-like-action-transitivas} Let $\G$ be an
			infinite, connected, and locally finite graph. Then $\G$ admits a
			transitive translation-like action by $\Z$ if and only if it has
			either one or two ends. Moreover, the action can be taken with $d(v,v\ast1)\leq3$
			for every vertex $v$. 
		\end{theorem}
		
		A problem left in \cite[Problem 3.5]{seward_burnside_2014} was to
		characterize which graphs admit a transitive translation-like action
		by $\Z$. Thus we have solved the case of locally finite graphs, and
		it only remains the case of graphs with vertices of infinite degree. 
		
		We now mention an application of these results to the problem of Hamiltonicity
		of Cayley graphs. This is related to a special case of Lovász conjecture
		which asserts the following: if $G$ is a finite group, then for every
		set of generators the associated Cayley graph admits a Hamiltonian
		path. Note that the existence of at least one such generating set
		is obvious ($S=G$), and the difficulty of the question, which is
		still open, is that it alludes every generating set. Now assume that
		$G$ is an infinite group, $S$ is a finite set of generators, and
		$\cay(G,S)$ admits a transitive translation-like action by $\Z$.
		This action becomes a bi-infinite Hamiltonian path after we enlarge
		the generating set, and thus it follows from Seward's theorem that
		every finitely generated group with one or two ends admits a generating
		set for which the associated Cayley graph admits a bi-infinite Hamiltonian
		path \cite[Theorem 1.8]{dumont_burnside_2018}. It is an open question
		whether this holds for every Cayley graph \cite[Problem 4.8]{dumont_burnside_2018},
		but our result yields an improvement in this direction. 
		
		\begin{corollary}
			Let $G$ be a finitely generated group with one or two ends, and let
			$S$ be a finite set of generators. Then the Cayley graph of $G$
			with respect to the generating set $\{g\in G|\ d_{S}(g,1_{G})\leq3\}$
			admits a bi-infinite Hamiltonian path. 
		\end{corollary}
		
		This was known to hold for generating sets of the form $\{g\in G|\ d_{S}(g,1_{G})\leq J\}$,
		where $S\subset G$ is a finite generating set for $G$ and $J$ depends
		linearly on the vertex degrees in $\cay(G,S)$. 
		
		In the more general case where we impose no restrictions on ends,
		we obtain the following result for non transitive translation-like
		actions. Observe that this readiliy implies \Cref{thm:seward}. 
		
		\begin{theorem}
			\label{thm:t2-translation-like-actions}Every infinite graph which
			is locally finite and connected admits a translation-like action by
			$\Z$. Moreover, the action can be taken with $d(v,v\ast1)\leq3$
			for every vertex $v$. 
		\end{theorem}
		
		These statements about translation-like actions can also be stated
		in terms of powers of graphs. Given a graph $\G$, its $n$-th power
		$\G^{n}$ is defined as the graph with the same set of vertices, and
		where two vertices $u,v$ are joined if their distance in $\G$ is
		at most $n$. It is well-known that the cube of every finite and connected
		graph is Hamiltonian \cite{zbMATH03278210,sekanina1960ordering,karaganis_cube_1968}.
		Our \Cref{thm:t1-translation-like-action-transitivas} generalizes
		this to infinite and locally finite graphs. That is, it shows that
		the cube of a locally finite and connected graph with one or two ends
		admits a bi-infinite Hamiltonian path. 
		
		We mention that \Cref{thm:t2-translation-like-actions} has been
		proved in \cite[Section 4]{cohen_strongly_2021}, using the same fact
		about cubes of finite graphs. 
		
		\subsection{Computability of translation-like actions}
		
		Now we turn our attention to the problem of computing translation-like
		actions on groups or graphs. We recall that a graph is computable
		if there exists an algorithm which given two vertices, determines
		whether they are adjacent or not. If moreover the graph is locally
		finite, and the function that maps a vertex to its degree is computable,
		then the graph is said to be highly computable. This extra condition
		is necessary to compute the neighborhood of a vertex.
		
		An important example comes from group theory: if $G$ is a finitely
		generated group with decidable word problem and $S$ is a finite set
		of generators, then its Cayley graph with respect to $S$ is highly
		computable.
		
		There is a variety of problems in graph theory that have no computable
		solution for infinite graphs. A classical example is the problem of
		computing infinite paths. Kőnig's infinity lemma asserts that every
		infinite, connected, and locally finite graph admits an infinite path.
		However, there are highly computable graphs which admit paths, all
		of which are uncomputable \cite{jockusch_pi_1972}. Another example
		comes from Hall's matching theorem. There are highly computable graphs
		satisfying the hypotheses in the theorem, but which admit no computable
		right perfect matching \cite{manaster_effective_1972}. These two
		results are used in the proof of Seward's theorem, so the translation-like
		actions from this proof are not clearly computable. We say that a
		translation-like action by $\Z$ on a graph is computable when there
		is an algorithm which given a vertex $v$ and $n\in\Z$, computes
		the vertex $v\ast n$. %

		Our interest in the computability of translation-like actions comes
		from symbolic dynamics, and the shift spaces associated to a group.
		We will need a computable translation-like action such that it is
		possible to distinguish in a computable manner between different orbits.
		We introduce here a general definition, though we will only treat
		the case where the acting group is $\Z$. 
		\begin{definition}
			\label{def:orbit-membership-problem}Let $G$ be a group, and let
			$S\subset G$ be a finite set of generators. A group action of $H$
			on $G$ is said to have \textit{decidable} \textit{orbit membership
				problem} if there exists an algorithm which given two words $u$ and
			$v$ in $(S\cup S^{-1})^{*}$, decides whether the corresponding group
			elements $u_{G},v_{G}$ lie in the same orbit under the action. 
		\end{definition}
		
		Note that if $H$ is a subgroup of $G$, then the action $H\curvearrowright G$
		by right translations has decidable orbit membership problem if and
		only if $H$ has decidable subgroup membership problem (\Cref{prop:obvio}).
		Thus this property can be regarded as the geometric reformulation,
		in the sense of Whyte \cite{whyte_amenability_1999}, of the subgroup
		property of having decidable membership problem. The orbit membership
		problem has been studied for some actions by conjugacy and by group
		automorphisms (see \cite{bogopolski_orbit_2009,burillo_conjugacy_2016,ventura_grouptheoretic_2014}
		and references therein).

		Our main result associated to computable translation-like actions
		on groups is the following. 
		
		\begin{theorem}
			\label{thm:computable-translation-like-actions-with-decidable-orbit-problem}Let
			$G$ be a finitely generated infinite group with decidable word problem.
			Then $G$ admits a translation-like action by $\Z$ that is computable
			and has decidable orbit membership problem.
		\end{theorem}
		
		The proof of \Cref{thm:computable-translation-like-actions-with-decidable-orbit-problem}
		proceeds as follows. For groups with one or two ends, we will show
		the existence of a computable and transitive translation-like action
		by $\Z$, that is, a computable version of \Cref{thm:t1-translation-like-action-transitivas}.
		This action has decidable orbit membership problem for the trivial
		reason that it has only one orbit. For groups with at least two ends
		we obtain the action from a subgroup. It follows from Stalling's structure
		theorem on ends of groups that a finitely generated groups with two
		or more ends has a subgroup isomorphic to $\Z$. We will show that,
		if the group has solvable word problem, then this subgroup has decidable
		membership problem. This proof is based on the computability of normal
		forms associated to Stalling's structure theorem (\Cref{prop:stallings-descomposicion-calculable}).

		\subsection{Medvedev degrees of effective subshifts}
		
		We now turn our attention to Medvedev degrees, a complexity measure
		which is defined using computable functions. Precise definitions of
		this and the following concepts are given in \Cref{sec:medvedev}.
		Intuitively, the Medvedev degree of a set $P\subset A^{\N}$ measures
		how hard is to compute a point in $P$. For example, a set has zero
		Medvedev degree if and only if it has a computable point. This complexity
		measure becomes meaningful when we regard $P$ as the set of solutions
		to a problem. This notion can be applied to a variety of objects,
		such as graph colorings \cite{remmel_graph_1986}, paths on graphs,
		matchings from Hall's matching theorem, and others \cite[Chapter 13]{ershov_handbook_1998a}.
		In this article we consider Medvedev degrees of subshifts. %
	
		Let $G$ be a finitely generated group, and let $A$ be a finite alphabet.
		A subshift is a subset of $A^{G}$ which is closed in the prodiscrete
		topology, and is invariant under translations. Dynamical properties
		of subshifts have been related to their computational properties in
		different ways. A remarkable example is the characterization of the
		entropies of two dimensional subshifts of finite type as the class
		of nonnegative $\Pi_{1}^{0}$ real numbers \cite{hochman_characterization_2010}.
		
		Here we adress the problem of classifying what Medvedev degrees can
		be attained for a certain class of subshifts. For instance, this classification
		is known for subshifts of finite type in the groups $\Z^{d}$, $d\geq1$.
		In the case $d=1$, all subshifts of finite type have Medvedev degree
		zero, because all of them contain a periodic point, and then a computable
		point. In the case $d\geq2$, subshifts of finite type can attain
		the class of $\Pi_{1}^{0}$ Medvedev degrees \cite{simpson_medvedev_2014}. 
		
		A larger class of subshifts is that of effective subshifts. A subshift
		over $\Z$ is effective if the set of words which do not appear in
		its configurations is computably enumerable. This notion can be extended
		to a finitely generated group, despite some intricacies that arise
		in relation to the word problem of the group. We will only deal with
		groups with decidable word problem, and the notion of effective subshift
		is a straightforward generalization. 
		
		Answering a question left open in \cite{simpson_medvedev_2014}, Joseph
		Miller proved that effective subshifts over $\Z$ can attain all $\Pi_{1}^{0}$
		Medvedev degrees \cite{miller_two_2012}. We generalize this result
		to the class of infinite, finitely generated groups with decidable
		word problem. 
		
		\begin{theorem}
			\label{thm:medvedev-degrees-of-effective-subshifts}Let $G$ be a
			finitely generated and infinite group with decidable word problem.
			The class of Medvedev degrees of effective subshift on $G$ is the
			class of all $\Pi_{1}^{0}$ Medvedev degrees.
		\end{theorem}
		
		The idea for the proof is the following. Given any subshift $Y\subset A^{\Z}$,
		we can construct a new subshift $X\subset B^{G}$ that simultaneously
		describes translation-like actions $\Z\curvearrowright G$, and elements
		in $Y$. Then \Cref{thm:computable-translation-like-actions-with-decidable-orbit-problem}
		ensures that this construction preserves the Medvedev degree of $Y$,
		and the result follows from the known classification for $\Z$ \cite{miller_two_2012}. 
		
		Despite the simplicity of the proof, we need to translate some computability
		notions from $A^{\N}$ to $A^{G}$, this is done by taking a computable
		numbering of $G$. The notions obtained do not depend of the numbering,
		and are compatible with notions already present in the literature
		that were defined by other means \cite{aubrun_notion_2017}. 
		
		This construction using translation-like actions was introduced in
		\cite{jeandel_translationlike_2015}, and has been used to prove different
		results in the context of symbolic dynamics. For example, to transfer
		results about the emptiness problem for subshifts of finite type from
		one group to another \cite{jeandel_translationlike_2015}, to produce
		aperiodic subshifts of finite type on new groups \cite{cohen_strongly_2021,jeandel_translationlike_2015},
		and to study the entropy of subshifts of finite type on some amenable
		groups \cite{barbieri_entropies_2021}.
		
		\subsection*{Paper structure}
		
		In \Cref{sec:Preliminaries} we fix some notation, and recall
		some basic facts on graph theory, group theory, and computability
		theory on countable sets. In \Cref{sec:Translation-like-actions-by}
		we show \Cref{thm:t1-translation-like-action-transitivas} and
		\Cref{thm:t2-translation-like-actions} about translation-like
		actions. In \Cref{sec:Computable-translation-like-acti} we show
		some results about computable translation-like actions, including
		\Cref{thm:computable-translation-like-actions-with-decidable-orbit-problem}.
		This result is applied in \Cref{sec:medvedev} to prove \Cref{thm:medvedev-degrees-of-effective-subshifts}
		on Medvedev degrees.
		
		\subsection*{Acknowledgements}
		
		This paper would not have been possible without the help, guidance,
		and careful reading of my two advisors. I am grateful to Sebastián
		Barbieri for suggesting that translation-like actions would work to
		prove \Cref{thm:medvedev-degrees-of-effective-subshifts}, and
		for helping me with many details. I am also grateful to Cristóbal
		Rojas for providing the computability background which took me to
		Medvedev degrees of complexity, and finally to ask \Cref{thm:medvedev-degrees-of-effective-subshifts}.
		Part of this work was done during a research stay at Institut de Mathématiques
		de Toulouse under the supervision of Mathieu Sablik, and I am also
		indebted to his hospitality. 
		
		This research was partially supported by ANID 21201185 doctorado nacional,
		ANID/Basal National Center for Artificial Intelligence CENIA FB210017,
		and the European Union's Horizon 2020 research and innovation program
		under the Marie Sklodowska-Curie grant agreement No 731143.
		
		\section{Preliminaries}\label{sec:Preliminaries}
		
		We denote by $f\circ g$ the function that applies $g$ to the argument,
		and then $f$.
		
		\subsection{Graph theory}
		
		In this article all graphs are undirected and unlabeled. Loops and
		multiple edges are allowed. The vertex set of a graph $\G$ will be
		denoted by $V(\G)$, and its edge set by $E(\G)$. Each edge \textit{joins
		}a pair of vertices, and is said to be \textit{incident }to them.
		Two vertices joined by an edge are called \textit{adjacent}. The \textit{degree}
		$\deg_{\G}(v)$ of the vertex $v$ is the number of incident edges
		to $v$, where loops are counted twice. A graph is said to be \textit{finite}
		when its edge set is finite, and \textit{locally finite} when every
		vertex has finite degree. 
		
		In our constructions we will constantly consider induced subgraphs.
		Given a set of vertices $V\subset V(\G)$, the \textit{induced subgraph
		}$\G[V]$ is the subgraph of $\G$ whose vertex set equals $V$, and
		whose edge set is that of all edges in $E(\G)$ whose incident vertices
		lie in $V$. On the other hand, $\G-V$ stands for the subgraph of
		$\G$ obtained by removing from $\G$ all vertices in $V$, and all
		edges incident to vertices in $V$. That is, $\G-V$ equals the induced
		subgraph $\G[V(\G)-V]$. If $\L$ is a subgraph of $\G$, we denote
		by $\G-\L$ the subgraph $\G-V(\L)$. 
		
		A \textit{path }on $\G$\textit{ }is an injective function $f\colon[a,b]\to V(\G)$
		that sends consecutive integers to adjacent vertices, where $[a,b]\subset\Z$.
		We introduce now some useful terminology for dealing with paths. We
		say that $f$ \textit{joins $f(a)$ }to $f(b)$, and define its \textit{length
		}as $b-a$. We say that $f$ \textit{visits }the vertices in its image,
		and we denote this set by $V(f)$. We denote by $\G-f$ the subgraph
		$\G-V(f)$. The vertices $f(a)$ and $f(b)$ are called the \textit{initial}
		and \textit{final} vertices of $f$, respectively. When every pair
		of vertices in the graph $\G$ can be joined by a path, then we say
		that $\G$ is \textit{connected}. In this case we define the \textit{distance}
		between two vertices as the length of the shortest path joining them.
		This distance induces the \textit{path-length }metric on $V(\G)$,
		which we denote by $d_{\G}$. 
		
		A \textit{connected component }of $\G$ is a connected subgraph of
		$\G$ which is maximal for the subgraph relation. The \textit{number
			of ends }of $\G$ is the supremum of the number of infinite connected
		components of $\G-V$, where $V$ ranges over all finite sets of vertices
		in $\G$.
		
		\subsection{Words and finitely generated groups}
		
		We now review some terminology and notation on words, alphabets, and
		finitely generated groups. An \textit{alphabet }is a finite set. The
		set of finite words on alphabet $A$ is denoted by $A^{*}$. The empty
		word is denoted by $\e$. A word $u$ of length $n$ is a \textit{prefix
		}of $v$ when they coincide in the first $n$ symbols. %
		
		Now let $G$ be a group. The identity element of $G$ is denoted by
		$1_{G}$, or $1$ if no confusion arises. Let $S\subset G$ be a finite
		set, and let $S^{-1}$ be the set of formal inverses to elements in
		$S$. Given a word $w\in(S\cup S^{-1})^{*}$, we denote by $w_{G}$
		the group element obtained by multiplying in $G$ the elements from
		$S$ that constitute the word. We also write $u=_{G}v$ when the words
		$u,v$ correspond to the same group element. A set $S\subset G$ is
		said to \textit{generate} $G$ if every group element can be written
		as a word in $(S\cup S^{-1})^{*}$, and $G$ is \textit{finitely generated}
		when it admits a finite generating set. A finite generating set $S\subset G$
		induces the \textit{left-invariant word metric }on $G$, denoted by
		$d_{S}$. The distance $d_{S}(g,h)$ is the length of the shortest
		word $\ensuremath{w\in(S\cup S^{-1})^{*}}$ such that $g(w)_{G}=h$.
		
		If $S\subset G$ is a finite set of generators, we denote by $\cay(G,S)$
		the (undirected, and right)\textit{ Cayley graph} of $G$ relative
		to $S$. The vertex set of $\cay(G,S)$ is $G$, and the edge set
		of $\cay(G,S)$ is $\{(g,gs)\mid g\in G,s\in S\cup S^{-1}\}$. The
		edge $(g,gs)$ joins the vertex $g$ with the vertex $gs$. Note that
		the distance that this graph assigns to a pair of elements in $G$
		equals their distance in the word metric associated to the same generating
		set. The \textit{number of ends }of a finitely generated group is
		the number of ends of its Cayley graph, for any generating set. This
		definition does not depend on the chosen generating set, and can only
		be among the numbers $\{0,1,2,\infty\}$ \cite{freudenthal_ueber_1945,hopf_enden_1944}. 
		
		We now recall some algorithmic properties of groups and subgroups.
		The concept of decidable set of words is defined in the next subsection.
		Let $S\subset G$ be a finite set of generators, and let $H$ be a
		subgroup of $G$. We say that $H$ has \textit{decidable subgroup
			membership problem }if $\{w\in(S\cup S^{-1})^{*}\mid w_{G}\in H\}$
		is a decidable subset of $(S\cup S^{-1})^{*}$. This notion does not
		depend on the chosen generating set. In the particular case where
		$H=\{1_{G}\}$, the set defined above is called the \textit{word problem}
		of $G$. The property of having decidable word problem is closely
		related to the property of being a computable group, which we discuss
		in more detail in the next subsection. 
		
		\subsection{Computability theory on countable sets via numberings}\label{subsec:numberings}
		
		We start by reviewing some classical notions from recursion theory
		or computability theory. All these facts are well-known, the reader
		is referred to \cite{chiswell_course_2009} for computability theory,
		and to \cite[Chapter 14]{griffor_handbook_1999} for a survey on numberings. 
		
		We will use the word \textit{algorithm} to refer to the formal object
		of Turing machine. We will use other common synonyms such as ``effective
		procedure''. A partial function $f\colon D\subset\N\to\N$ is \textit{computable
		}if there is an algorithm satisfying the following. On input $n$,
		the algorithm halts if and only if $n\in D$, and in this case outputs
		$f(n)$. A subset $D\subset\N$ is \textit{semi-decidable }when there
		is an algorithm that halts on input $n$ if and only if $n\in D$.
		A set $D\subset\N$ is \textit{decidable }when both $D$ and $\N-D$
		are semi-decidable.
		
		All these notions extend directly to products $\N^{p}$, $p\geq1$,
		and sets of words $A^{\ast}$, as these objects can be represented
		by natural numbers in a canonical way. In order to extend these notions
		to other objects such as graphs and countable groups, we take a unified
		approach via numberings:
		\begin{definition}
			A (bijective)\textit{ numbering }of a set $X$ is a bijective map
			$\nu:N\to X$, where $N$ is a decidable subset of $\N$. We call
			$(X,\nu)$ a \textit{numbered set}. When $\nu(n)=x$, we say that
			$n$ is a \textit{name} for $x$, or that $n$ represents $x$. 
		\end{definition}
		
		A numbering of $X$ defines computability notions in $X$ in the same
		manner that charts are used to define continuous or differentiable
		function on manifolds. For instance, a function $f\colon X\to X$
		is computable on $(X,\nu)$ when the ``function in charts'' $\nu^{-1}\circ f\circ\nu$
		is computable. There is a notion of equivalence for numberings: two
		numberings $\nu$, $\nu'$ of $X$ are equivalent when the identity
		function $(X,\nu)\to(X,\nu')$ is computable. The Cartesian product
		$X\times X'$ of two numbered sets $(X,\nu)$, $(X',\nu')$ admits
		a unique numbering -up to equivalence- for which the projection functions
		to $(X,\nu)$, $(X',\nu')$ are computable. This provides definitions
		of computable functions and relations between different numbered sets,
		and we can freely speak about computable functions and relations between
		numbered sets. We will be interested in the following objects:
		
		\begin{definition}
			A graph $\G$ is computable if we can endow $V(\G)$ and $E(\G)$
			with numberings, in such a manner that the relation of adjacency,
			and the relation $\{(e,u,v)\mid e\text{ joins }u\text{ and }v\}$
			are decidable. We say that $\G$ is also \textit{highly computable}
			when it is locally finite, and the vertex degree function $V(\G)\to\N$,
			$v\mapsto\deg_{\G}(v)$ is computable. 
		\end{definition}
		
		\begin{definition}
			A numbering $\nu$ of a group $G$ is said to be \textit{computable}
			when it makes the group operation $G^{2}\to G$ is computable. In
			this case, the pair $(G,\nu)$ is called a \textit{computable group. }
		\end{definition}
		
		These notions provide a formal and precise meaning to general statements
		about the computability of objects such as translation-like actions
		and bi-infinite paths on computable groups or graphs. For instance,
		a group action on a computable group $G$ is computable when the function
		$(g,n)\to g\ast n$ from the numbered set $G\times\Z$ to the numbered
		set $G$ is computable.
		
		It is well known that algorithmic properties of finitely generated
		groups have a number of stability properties, such as being independent
		of the generating set. In terms of numberings, this is expressed as
		follows:
		\begin{proposition}
			\label{prop:computabilitiy-groups}Let $G$ be a finitely generated
			group. Then:
			\begin{enumerate}
				\item $G$ admits a computable numbering if and only if it has decidable
				word problem.
				\item If $G$ admits a computable numbering, then all computable numberings
				of $G$ are equivalent.
				\item If $H$ is another finitely generated computable group, then any group
				homomorphism $f\colon G\to H$ is computable. 
			\end{enumerate}
		\end{proposition}
		
		\begin{proof}[Proof sketch]
			Suppose that $G$ has decidable word problem, let $S\subset G$
			be a finite generating set, and let $\pi\colon(S\cup S^{-1})^{\ast}\to G$
			be the function that sends a word to the corresponding group element.
			Using the decidability of the word problem, we can compute a set $N\subset(S\cup S^{-1})^{\ast}$
			such that the restriction of $\pi$ to $N$ is a bijection. Being
			$N$ a decidable subset of $(S\cup S^{-1})^{\ast}$, it admits a computable
			bijection with $\N$. The composition of these functions give a bijection
			$\nu\colon\N\to G$, and it is easy to verify that it is a computable
			numbering. The reverse implication is left to the reader. Items $2$
			and $3$ are also left to the reader: the relevant functions are determined
			by the finite information of letter-to-word substitutions, and this
			allows to prove that they are computable. 
		\end{proof}
		We will also make use of the following well-known fact. The proof
		is straightforward, and left to the reader. 
		\begin{proposition}
			\label{prop:cayley-graph-of-computable-group-is-computable} Let $G$
			be a finitely generated group with decidable word problem, and let
			$S$ be a finite generating set. Then $\cay(G,S)$ is a highly computable
			graph. 
		\end{proposition}

		\section{Translation-like actions by $\Z$ on locally finite graphs}\label{sec:Translation-like-actions-by}
		
		The goal of this section is to prove \Cref{thm:t1-translation-like-action-transitivas}
		and \Cref{thm:t2-translation-like-actions}. That is, that every
		connected, locally finite, and infinite graph admits a translation
		by $\Z$, and that this action can be chosen transitive exactly when
		the graph has one or two ends. The actions that we construct satisfy
		that the distance between a vertex $v$ and $v\ast1$ is at most $3$. 
		
		Our proof goes by constructing these actions locally, and in terms
		of 3-paths: 
		\begin{definition}
			Let $\G$ be a graph. A \textit{3-path} on $\G$ is an injective function
			$f\colon[a,b]\to V(\G)$ such that consecutive integers in $[a,b]$
			are mapped to vertices whose distance is at most 3. A \textit{bi-infinite
			}3-path on $\G$ is an injective function $f\colon\Z\to V(\G)$ satisfying
			the same condition on the vertices. A 3-path or bi-infinite 3-path
			is called \textit{Hamiltonian }when it is also a surjective function. 
		\end{definition}
		
		It is well known that every finite and connected graph admits a Hamiltonian
		3-path, where we can choose its initial and final vertex \cite{zbMATH03278210,sekanina1960ordering,karaganis_cube_1968}.
		Here we will need a slight refinement of this fact:

			\begin{lemma}\label{lem:lema-tecnico-finite-3-paths} Let $\G$ be a graph that
			is connected and finite. For every pair of different vertices $u$
			and $v$, $\G$ admits a Hamiltonian $3$-path $f$ which starts at
			$u$, ends at $v$, and moreover satisfies the following two conditions:
			\begin{enumerate}
				\item The first and last ``jump'' have length at most $2$. That is, if
				$f$ visits $w$ immediately after the initial vertex $u$, then $d_{\G}(u,w)\leq2$.
				Moreover, if $f$ visits $w$ immediately before the final vertex
				$v$, then $d_{\G}(w,v)\leq2$. , 
				\item There are no consecutive ``jumps'' of length $3$. That is, if $f$
				visits $w_{1}$, $w_{2}$ and $w_{3}$ consecutively, then $d_{\G}(w_{1},w_{2})\leq2$
				or $d_{\G}(w_{2},w_{3})\leq2$. 
			\end{enumerate}
		\end{lemma}
		
		Let us review some terminology on 3-paths before proving this result.
		When dealing with 3-paths, we will use the same terms introduced for
		paths in the preliminaries, such as initial vertex, final vertex,
		visited vertex, etc. Let $f$ and $g$ be 3-paths. We say that $f$
		extends $g$ if its restriction to the domain of $g$ equals $g$.
		We will extend 3-paths by concatenation, which we define as follows.
		Suppose that final vertex of $f$ is at distance at most 3 from the
		initial vertex of $g$, and such that $V(f)\cap V(g)=\emptyset$.
		The \textit{concatenation of $f$, $g$ }is the 3-path that extends
		$f$, and after the final vertex of $f$ visits all vertices visited
		by $g$ in the same order. Finally, the \textit{inverse }of the 3-path
		$f$, denoted by $-f$, is defined by $(-f)(n)=f(-n)$. Note that
		its domain is also determined by this expression. 
		\begin{proof}[Proof of \Cref{lem:lema-tecnico-finite-3-paths}]
The proof is by induction of the cardinality of $V(\G)$. The claim clearly holds if $|V(\G)|\leq 2$. Now assume that $\G$ is a connected finite graph with $|V(\G)|\geq 3$, and let $u$ and $v$ be two different vertices. We consider the connected components of the graph $\G-\{v\}$ obtained by removing the vertex $v$ from $\G$. Let $\G_{u}$ be the finite connected component of $\G-\{v\}$ that contains $u$, and let $\G_{v}$ be the subgraph of $\G$ induced by the set of vertices $V(\G)-V(\G_{u})$. Thus $u\in\G_{u}$, $v\in\G_{v}$, and both $\G_{u}$ and $\G_{v}$ are connected. Let us first assume that both $\G_{u}$ and $\G_{v}$
are graphs with at least two vertices. Then we can apply the inductive hypothesis on each one of them. Let $f_{u}$ be a Hamiltonian 3-path on $\G_{u}$ as in the statement, whose initial vertex is $u$, and whose final vertex $u'$ is at distance to $v$ at most 2. Let $f_{v}$ be a Hamiltonian 3-path on $\G_{v}$ as in the statement, whose initial vertex $v'$ is adjacent to $v$, and whose final vertex is $v$.  We claim that the 3-path $f$ obtained by concatenating $f_{u}$, $f_{v}$ verifies the required conditions. It is clear that  $d_{\Gamma}(u',v')\leq 3$, and thus $f$ is a 3-path. It is also clear that $f$ verifies the first condition in the statement. Regarding the second condition, it suffices to show that the ``jump'' from $u'$ to $v'$ is between two ``jumps'' with length at most two. That is, that the vertex visited by $f$ before $u'$ is at distance at most 2 from $u'$, and that the vertex visited by $f$ after $v'$ is at distance at most $2$ from $v'$. Indeed, this follows from the fact that both $f_u$ and $f_v$ verify the first condition in the statement. This finishes the argument in the case that both $\G_u$ and $\G_v$ have at least two vertices. If $\G_u$ has one vertex and $\G_v$ has at least two vertices then we modify the previous procedure by redefining $f_u$ as the 3-path that only visits $u$. It is easy to verify that then the concatenation of $f_u$ and $f_v$ is a 3-path in $\G$ verifying the two numbered conditions. The case where $\G_v$ has one vertex and $\G_u$ has at least two vertices is symmetric, and the case where both $\G_u$ and $\G_v$ have one vertex is excluded since we assumed $|V(G)|\geq 3$.
		\end{proof}
		We will define bi-infinite 3-paths by extending finite ones iteratively.
		The following definition will be key for this purpose:
		\begin{definition}
			\label{def:bi-right-extensible}Let $f$ be a 3-path on a graph $\G$.
			We say that $f$ is \textit{bi-extensible} if the following conditions
			are satisfied:
			\begin{enumerate}
				\item $\G-f$ has no finite connected component.
				\item There is a vertex $u$ in $\G-f$ at distance at most 3 from the final
				vertex of $f$. 
				\item There is a vertex $v\ne u$ in $\G-f$ at distance at most 3 from
				the initial vertex of $f$. 
			\end{enumerate}
			If only the first two conditions are satisfied, we say that $f$ is
			\textit{right-extensible}.
		\end{definition}
		
		We will now prove some elementary facts about the existence of 3-paths
		that are bi-extensible and right-extensible. The proofs are elementary,
		and are given by completeness.
		\begin{lemma}
			\label{lem:existencia-right-extensible}Let $\G$ be a graph that
			is infinite, connected, and locally finite. Then for every pair of
			vertices $u$ and $v$ in $\G$, there is a right-extensible 3-path
			whose initial vertex is $u$, and which visits $v$. 
		\end{lemma}
		
		\begin{proof}
			As $\G$ is connected, there is a path $f$ joining $u$ to $v$.
			Now define $\L$ as the graph induced in $\G$ by the set of vertices
			that are visited by $f$, or that lie in a finite connected component
			of $\G-f$. Notice that as $\G$ is locally finite, there are finitely
			many such connected components, and thus $\L$ is a finite and connected
			graph. By construction, $\G-\L$ has no finite connected component. 
			
			The desired 3-path will be obtained as a Hamiltonian 3-path on $\L$.
			Indeed, as $\G$ is connected, there is a vertex $w$ in $\L$ that
			is adjacent to some vertex in $\G-\L$. By \Cref{lem:lema-tecnico-finite-3-paths}
			there is a 3-path $f'$ which is Hamiltonian on $\L$, starts at $u$
			and ends in $w$. We claim that $f'$ is right-extensible. Indeed,
			our choice of $\L$ ensures that $\G-f'$ has no finite connected
			component, and our choice of $w$ ensures that the final vertex of
			$f'$ is adjacent to a vertex in $\G-f'$. 
		\end{proof}
		\begin{lemma}
			\label{lem:existencia-bi-extensible} Let $\G$ be a graph that is
			infinite, connected, and locally finite. Then for every vertex $u$
			in $\G$, there is a bi-extensible 3-path in $\G$ that visits $u$. 
		\end{lemma}
		
		\begin{proof}
			Let $v$ be a vertex in $\G$ that is adjacent to $u$, with $v\neq u$.
			Let $\L$ be the subgraph of $\G$ induced by the set of vertices
			that lie in a finite connected component of $\G-\{u,v\}$, or in $\{u,v\}$.
			As $\G$ is locally finite, there are finitely many such connected
			components, and thus $\L$ is a finite and connected subgraph of $\G$.
			By construction, $\G-\L$ has no finite connected component. 
			
			The desired 3-path will be obtained as a Hamiltonian 3-path on $\L$.
			Indeed, as $\G$ is connected there are two vertices $w\in V(\G-\L)$
			and $w'\in V(\L)$, with $w$ adjacent to $w'$ in $\G$. As $\L$
			has at least two vertices, we can invoke \Cref{lem:lema-tecnico-finite-3-paths}
			to obtain a 3-path $f$ that is Hamiltonian on $\L$, whose initial
			vertex is $w'$, and whose final vertex is adjacent to $w'$. It is
			clear that then $f$ is a bi-extensible 3-path in $\G$. 
		\end{proof}
		Our main tool to construct bi-infinite 3-paths is the following result,
		which shows that bi-extensible 3-paths can be extended to larger bi-extensible
		3-paths.
		\begin{lemma}
			\label{lem:extend-3-paths}Let $\G$ be a graph that is infinite,
			connected, and locally finite. Let $f$ be a bi-extensible 3-path
			on $\G$, and let $u$ and $v$ be two different vertices in $\G-f$
			whose distance to the initial and final vertex of $f$ is at most
			3, respectively. Let $w$ be a vertex in the same connected component
			of $\G-f$ that some of $u$ or $v$. Then there is a 3-path $f'$
			which extends $f$, is bi-extensible on $\G$, and visits $w$. Moreover,
			we can assume that the domain of $f'$ extends that of $f$ in both
			directions.
		\end{lemma}
		
		\begin{proof}
			If $u$ and $v$ lie in different connected components of $\G-f$,
			then then the claim is easily obtained by applying \Cref{lem:existencia-right-extensible}
			on each of these components. Indeed, by \Cref{lem:existencia-right-extensible}
			there are two right-extensible 3-paths $g$ and $h$ in the corresponding
			connected components of $\G-f$, such that the initial vertex of $g$
			is $u$, the initial vertex of $h$ is $v$, and some of them visits
			$w$. Then the concatenation of $-g$, $f$ and $h$ satisfies the
			desired conditions. 
			
			We now consider the case where $u$ and $v$ lie in the same connected
			component of $\G-f$. This graph will be denoted $\L$. Note that
			$\L$ is infinite because $f$ is bi-extensible. We claim that there
			are two right-extensible 3-paths on $\L$, $g$ and $h$, satisfying
			the following list of conditions: the initial vertex of $g$ is $u$,
			the initial vertex of $h$ is $v$, some of them visits $w$, and
			$V(g)\cap V(h)=\emptyset$. In addition, $(\L-g)-h$ has no finite
			connected component, and has two different vertices $u'$ and $v'$
			such that $u'$ is at distance at most 3 from the last vertex of $g$,
			and $v'$ is at distance at most 3 from the last vertex of $h$. Suppose
			that we have $g,h$ as before. Then we can define a 3-path $f'$ by
			concatenating $-g$, $f$ and then $h$. It is clear that then $f'$
			satisfies the conditions in the statement. 
			
			We now construct $g$ and $h$. We start by taking a connected finite
			subgraph $\L_{0}$ of $\L$ which contains $u,v,w$ and such that
			$\L-\L_{0}$ has no finite connected component. The graph $\L_{0}$
			can be obtained, for instance, as follows. As $\Lambda$ is connected,
			we can take a path $f_{u}$ from $u$ to $w$, and a path $f_{v}$
			from $v$ to $w$. Then define $\L_{0}$ as the graph induced by the
			vertices in $V(f_{v})$, $V(f_{u})$, and all vertices in the finite
			connected components of $(\L-f_{v})-f_{u}$. 
			
			Let $p$ be a Hamiltonian 3-path on $\L_{0}$ from $u$ to $v$, as
			in \Cref{lem:lema-tecnico-finite-3-paths}. The desired 3-paths
			$f$ and $g$ will be obtained by ``splitting'' $p$ in two. As
			$\L$ is connected, there are two vertices $u_{0}\in V(\L_{0})$,
			$v'\in V(\L-\L_{0})$ such that $u_{0}$ and $v$ are adjacent in
			$\L$. By the conditions in \Cref{lem:lema-tecnico-finite-3-paths},
			there is a vertex $v_{0}$ in $V(\L_{0})$ whose distance from $u_{0}$
			is at most 2, and such that $p$ visits consecutively $\{u_{0},v_{0}\}$.
			We will assume that $p$ visits $v_{0}$ after visiting $u_{0}$,
			the other case being symmetric. As $\L-\L_{0}$ has no finite connected
			component, there is a vertex $u'$ in $\L-\L_{0}$ that is adjacent
			to $v'$. Thus, $u_{0}$ is at distance at most 2 from $u'$, and
			$v_{0}$ is at distance at most 3 from $v'$. Now we define $g$ and
			$h$ by splitting $p$ at the vertex $u_{0}$. More precisely, let
			$[a,c]$ be the domain of $p$, and let $b$ be such that $p(b)=u_{0}$.
			Then $h$ is defined as the restriction of $p$ to $[a,b]$, and we
			define $g$ by requiring $-g$ to be the restriction of $p$ to $[b+1,c].$
			Thus $h$ is a 3-path from $v$ to $v_{0}$, and $g$ is a 3-path
			from $u$ to $u_{0}$. By our choice of $\L_{0}$ and $p$, the 3-paths
			$h$ and $g$ satisfy the mentioned list of conditions, and thus the
			proof is finished. 
		\end{proof}
		When the graph has one or two ends, the hypotheses of \Cref{lem:extend-3-paths}
		on $u,v$ and $w$ are trivially satisfied. We obtain a very simple
		and convenient statement: we can extend a bi-extensible 3-path so
		that it visits a vertex of our choice. 
		\begin{corollary}
			\label{cor:one-two-ends}Let $\G$ be a graph that is infinite, connected,
			locally finite, and whose number of ends is either $1$ or $2$. Then
			for every bi-extensible 3-path $f$ and vertex $w$, there is a bi-extensible
			3-path on $\G$ that extends $f$ and visits $w$. We can assume that
			the domain of the new 3-path extends that of $f$ in both directions. 
		\end{corollary}
		
		We are now in position to prove some results about bi-infinite 3-paths.
		We start with the Hamiltonian case, which is obtained by iterating
		\Cref{cor:one-two-ends}. When we deal with bi-infinite 3-paths,
		we use the same notation and abbreviations introduced before for 3-paths,
		as long as they are well defined. 
		\begin{proposition}
			\label{prop:Hamiltonian-infinite-3paths}Let $\G$ be a graph that
			is infinite, connected, locally finite, and whose number of ends is
			either $1$ or $2$. Then $\G$ admits a bi-infinite Hamiltonian 3-path. 
		\end{proposition}
		
		\begin{proof}
			Let $(v_{n})_{n\in\N}$ be a numbering of the vertex set of $\G$.
			We define a sequence of bi-extensible 3-paths $(f_{n})_{n\in\N}$
			on $\G$ recursively. We define $f_{0}$ as a bi-extensible 3-path
			which visits $v_{0}$. The existence of $f_{0}$ is guaranteed by
			\Cref{lem:existencia-bi-extensible}. Now let $n\geq0$, and
			assume that we have defined a 3-path $f_{n}$ that visits $v_{n}$.
			We define $f_{n+1}$ as a bi-extensible 3-path on $\G$ which extends
			$f_{n}$, visits $v_{n+1}$, and whose domain extends the domain of
			$f_{n}$ in both directions. The existence such a 3-path is guaranteed
			by \Cref{cor:one-two-ends}. We have obtained a sequence $(f_{n})_{n\in\N}$
			such that for all $n$, $f_{n}$ visits $v_{n}$, and $f_{n+1}$ extends
			$f_{n}$. With this sequence we define a bi-infinite 3-path $f\colon\Z\to V(\G)$
			by setting $f(k)=f_{n}(k)$, for $n$ big enough. Note that $f$ is
			well defined because $f_{n+1}$ extends $f_{n}$ as a function, and
			the domains of $f_{n}$ exhaust $\Z$. By construction, $f$ visits
			every vertex exactly once, and thus it is Hamiltonian. 
		\end{proof}
		We now proceed with the non Hamiltonian case, where there are no restrictions
		on ends. We first prove that we can take a bi-infinite 3-path whose
		deletion leaves no finite connected component. 
		
		\begin{lemma}
			\label{lem:bi-infinte-3-path-removable}Let $\G$ be a graph that
			is infinite, connected, and locally finite. Then for every vertex
			$v$, there is a bi-infinite 3-path $f$ that visits $v$, and such
			that $\G-f$ has no finite connected component. 
		\end{lemma}
		
		\begin{proof}
			By \Cref{lem:existencia-bi-extensible} and \Cref{lem:extend-3-paths},
			$\G$ admits a sequence $(f_{n})_{n\in\N}$ of bi-extensible 3-paths
			such that $f_{0}$ visits $v$, $f_{n+1}$ extends $f_{n}$ for all
			$n\geq0$, and such that their domains exhaust $\Z$. We define a
			bi-infinite 3-path $f\colon\Z\to\G$ by setting $f(k)=f_{n}(k)$,
			for $n$ big enough. We claim that $f$ satisfies the condition in
			the statement, that is, that $\G-f$ has no finite connected component.
			We argue by contradiction. Suppose that $\G_{0}$ is a nonempty and
			finite connected component of $\G-f$. Define $V_{1}$ as the set
			of vertices in $\G$ that are adjacent to some vertex in $\G_{0}$,
			but which are not in $\G_{0}$. Then $V_{1}$ is nonempty as otherwise
			$\G$ would not be connected, and it is finite because $\G$ is locally
			finite. Moreover, $f$ visits all vertices in $V_{1}$, for otherwise
			$\G_{0}$ would not be a connected component of $\G-f$. As $V_{1}$
			is finite, there is a natural number $n_{1}$ such that $f_{n_{1}}$
			has visited all vertices in $V_{1}$. By our choice of $V_{1}$ and
			$n_{1}$, $\G_{0}$ is a nonempty and finite connected component of
			$\G-f_{n_{1}}$, and this contradicts the fact that $f_{n_{1}}$ is
			bi-extensible.
		\end{proof}
		Now the proof of the following result is by iteration of \Cref{lem:bi-infinte-3-path-removable}. 
		\begin{proposition}
			\label{prop:existencia-3-paths}Let $\G$ be a graph that is infinite,
			connected, and locally finite. Then there is a collection of bi-infinite
			3-paths $f_{i}\colon\Z\to\G$, $i\in I$, such that $V(\G)$ is the
			disjoint union of $V(f_{i})$, $i\in I$. 
		\end{proposition}
		
		\begin{proof}
			By \Cref{lem:bi-infinte-3-path-removable}, $\G$ admits a bi-infinite
			3-path $f_{0}$ such that $\G-f_{0}$ has no finite connected component.
			Each connected component of $\G-f_{0}$ is infinite, and satisfies
			the hypotheses of \Cref{lem:bi-infinte-3-path-removable}. Thus
			we can apply \Cref{lem:bi-infinte-3-path-removable} on each
			of these connected components. Iterating this process in a tree-like
			manner, we obtain a family of 3-paths $f_{i}\colon\Z\to\G$, $i\in I$
			whose vertex sets $V(f_{i})$ are disjoint. As \Cref{lem:bi-infinte-3-path-removable}
			allows us to choose a vertex to be visited by the bi-infinite 3-path,
			we can choose $f_{i}$ ensuring that every vertex of $\G$ is visited
			by some $f_{i}$. In this manner, $V(\G)$ is the disjoint union of
			$V(f_{i})$, ranging $i\in I$.
		\end{proof}
		Finally, we can obtain \Cref{thm:t1-translation-like-action-transitivas}
		and \Cref{thm:t2-translation-like-actions} from our statements
		in terms of bi-infinite 3-paths.
		
		\begin{proof}[Proof of \Cref{thm:t1-translation-like-action-transitivas}]
			Let $\G$ be a graph as in the statement. By \Cref{prop:Hamiltonian-infinite-3paths},
			$\G$ admits a Hamiltonian bi-infinite 3-path $f$. We define a translation-like
			$\ast\colon V(\G)\times\Z\to V(\G)$ by the expression $v\ast n=f(f^{-1}(v)+n),\ n\in\Z$.
			This translation-like action is transitive because $f$ is Hamiltonian,
			and satisfies $d_{\G}(v,v\ast1)\leq3$ because $f$ is a bi-infinite
			3-path. 
			
			We now prove the remaining implication of the result. That is, that
			a connected and locally finite graph which admits a transitive translation-like
			action by $\Z$ must have either one or two ends. This is stated in
			\cite[Theorem 3.3]{seward_burnside_2014} for graphs with uniformly
			bounded vertex degree, but the same proof can be applied to locally
			finite graphs. For completeness, we provide an alternative argument.
			Let $\G$ be a connected and locally finite graph which admits a transitive
			translation-like action by $\Z$, denoted $\ast$. As the action is
			free, $V(\G)$ must be infinite, and thus $\G$ has at least $1$
			end. Suppose now that it has at least 3 ends to obtain a contradiction.
			Let $J=\max\{d_{\G}(v,v\ast1)\mid v\in V(\G)\}$. As $\G$ has at
			least $3$ ends, there is a finite set of vertices $V_{0}$ such that
			$\G-V_{0}$ has at least three infinite connected components, which
			we denote by $\G_{1},$ $\G_{2}$ and $\G_{3}$. By enlarging $V_{0}$
			if necessary, we can assume that any pair of vertices $u$ and $v$
			that lie in different connected components in $\G-V_{0}$, are at
			distance $d_{\G}$ at least $J+1$. Now as $V_{0}$ is finite, there
			are two integers $n\leq m$ such that $V_{0}$ is contained in $\{v\ast k\mid n\leq k\leq m\}$.
			By our choice of $V_{0}$, it follows that the set $\{v\ast k\mid k\geq m+1\}$
			is completely contained in one of $\G_{1}$, $\G_{2}$, or $\G_{3}$.
			The same holds for $\{v\ast k\mid k\leq n-1\}$, and thus one of $\G_{1}$,
			$\G_{2}$, or $\G_{3}$ must be empty, a contradiction. 
		\end{proof}
		
		\begin{proof}[Proof of \Cref{thm:t2-translation-like-actions}]
			Let $\G$ be a graph as in the statement, and let $f_{i},\ i\in I$ as in
			\Cref{prop:existencia-3-paths}. We define $\ast:V(\G)\times\Z\to V(\G)$
			by the expression
			\begin{align*}
				v\ast n & =f(f^{-1}(v)+n),\ n\in\Z,
			\end{align*}
			where $f$ is the only $f_{i}$ such that $v$ is visited by $f_{i}$.
			Observe that $v\ast1$ is well defined because $V(\G)=\bigsqcup_{i\in I}V(f_{i})$.
			This defines a translation-like action by $\Z$, where the distance
			from $v$ to $v\ast1$, $v\in V(\G)$ is uniformly bounded by 3.
		\end{proof}
		\begin{remark}
			The proof given in this section is closely related to the characterization
			of those infinite graphs that admit infinite Eulerian paths. This
			is a theorem of Erdős, Grünwald, and Weiszfeld \cite{erdos_eulerian_1936}.
			In the recent work \cite{carrasco-vargas_infinite_2024}, the author
			of this article gave a different proof of the Erdős, Grünwald, and
			Weiszfeld theorem, that complements the original result by also characterizing
			those finite paths that can be extended to infinite Eulerian ones.
			This characterization is very similar to the notion of bi-extensible
			defined here. Indeed, the proofs of \Cref{prop:Hamiltonian-infinite-3paths}
			and the proof of the mentioned result about Eulerian paths follow
			the same iterative construction. 
		\end{remark}
		
		\begin{remark}
			As we mentioned before, it is known that the cube of every finite
			and connected graph is Hamiltonian \cite{zbMATH03278210,sekanina1960ordering,karaganis_cube_1968}.
			\Cref{prop:Hamiltonian-infinite-3paths} can be considered as
			a generalization of this fact to locally finite graphs. That is, \Cref{prop:Hamiltonian-infinite-3paths}
			shows that the cube of every locally finite and connected graph with
			either $1$ or $2$ ends admits a bi-infinite Hamiltonian path. 
		\end{remark}
		
		We end this section by rephrasing a problem left in \cite[Problem 3.5]{seward_burnside_2014}.
		\begin{problem}
			Find necessary and sufficient conditions for a connected graph to
			admit a transitive translation-like action by $\Z$. 
		\end{problem}
		
		We have shown that for locally finite graphs, the answer to this problem
		is as simple as possible, involving only the number of ends of the
		graph. The problem is now open for graphs that are not locally finite.
		We observe that beyond locally finite graphs there are different and
		non-equivalent notions of ends \cite{diestel_graphtheoretical_2003},
		and thus answering the problem above also requires to determine which
		is the appropriate notion of ends. 
		
		\section{Computable translation-like actions by $\Z$}\label{sec:Computable-translation-like-acti}
		
		The goal of this section is to prove \Cref{thm:computable-translation-like-actions-with-decidable-orbit-problem}.
		Namely, that every finitely generated infinite group with decidable
		word problem admits a translation-like action by $\Z$, with the additional
		property of being computable and with decidable orbit membership problem. 
		
		The proof of \Cref{thm:computable-translation-like-actions-with-decidable-orbit-problem}
		is as follows. For groups with at most two ends, we prove the existence
		of a computable and transitive translation-like action. That is, we
		prove a computable version of \Cref{thm:t1-translation-like-action-transitivas}.
		For groups with more than two ends, we prove the existence of a subgroup
		isomorphic to $\Z$ and with decidable subgroup membership problem.
		Thus for groups with two ends we provide two different proofs for
		\Cref{thm:computable-translation-like-actions-with-decidable-orbit-problem}.
		A group with two ends is virtually $\Z$, and it would be easy to
		give a direct proof, but the intermediate statements may have independent
		interest (\Cref{thm:computable-t1} and \Cref{prop:stallings-descomposicion-calculable}). 
		
		\subsection{Computable and transitive translation-like actions by $\Z$}
		
		The goal of this subsection to prove that \Cref{thm:t1-translation-like-action-transitivas}
		is computable on highly computable graphs:
		
		\begin{theorem}[Computable \Cref{thm:t1-translation-like-action-transitivas}]
			\label{thm:computable-t1} Let $\G$ be a graph that is highly computable,
			connected, and has either $1$ or $2$ ends. Then $\G$ admits a computable
			and transitive translation-like action by $\Z$, where the distance
			between a vertex $v$ and $v\ast1$ is uniformly bounded by 3.
		\end{theorem}
		
		We start by proving that the bi-extensible property (\Cref{def:bi-right-extensible})
		is algorithmically decidable on highly computable graphs with one
		end.
		\begin{proposition}
			\label{prop:bi-extensible-decidible-1-end}Let $\G$ be a graph that
			is highly computable, connected, and has one end. Then it is algorithmically
			decidable whether a 3-path $f$ is bi-extensible. 
		\end{proposition}
		
		\begin{proof}
			It is clear that the second and third conditions in the definition
			of bi-extensible are algorithmically decidable. For the first condition,
			note that as $\G$ has one end, we can equivalently check whether
			$\G-f$ is connected. This is proved to be a decidable problem in
			\cite{carrasco-vargas_infinite_2024}, Lemma 5.6. Note that the mentioned
			result concerns the remotion of edges instead of vertices, but indeed
			this is stronger: given $f$, we compute the set $E$ of all edges
			incident to a vertex in $V(f)$, and then use \cite[Lemma 5.6]{carrasco-vargas_infinite_2024}
			with input $E$.
		\end{proof}
		For graphs with two ends we prove a similar result, but we need an
		extra assumption. 
		\begin{proposition}
			\label{prop:bi-extensible-decidible-2-end}Let $\G$ be a graph that
			is highly computable, connected and has two ends. Let $f_{0}$ be
			a bi-extensible 3-path on $\G$, such that $\G-f_{0}$ has two infinite
			connected components. Then there is an algorithm that on input a 3-path
			$f$ that extends $f_{0}$, decides whether $f$ is bi-extensible. 
		\end{proposition}
		
		\begin{proof}
			It is clear that the second and third conditions in the definition
			of bi-extensible are algorithmically decidable. We address the first
			condition. We prove the existence of a procedure that, given a 3-path
			$f$ as in the statement, decides whether $\G-f$ has no finite connected
			component. Given $f$, we start by computing the set $E$. In \cite[Lemma 5.5]{carrasco-vargas_infinite_2024}
			there is an effective procedure that halts if and only if $\G-f$
			has some finite connected component (the mentioned result mentions
			the remotion of edges instead of vertices, but indeed this is stronger:
			given $f$, we compute the set $E$ of all edges incident to a vertex
			in $V(f)$, and then use \cite[Lemma 5.5]{carrasco-vargas_infinite_2024}
			with input $E$).
			
			Thus we need an effective procedure that halts if and only if $\G-f$
			has no finite connected component. As $f$ extends $f_{0}$, this
			is equivalent to ask whether $\G-f$ has at most two connected components.
			The procedure is as follows: given $f$, we start by computing the
			set $V_{0}$ of vertices in $\G-f$ that are adjacent to a vertex
			visited by $f$. Then for every pair of vertices in $u,v\in V_{0}$,
			we search exhaustively for a path that that joins them, and that never
			visits vertices in $V(f)$. That is, a path in $\G-f$. Such a path
			will be found if and only if the connected component of $\G-f$ that
			contains $u$ equals the one that contain $v.$ We stop the procedure
			once we have found enough paths to write $V_{0}$ as the disjoint
			union $V_{1}\sqcup V_{2}$, where every pair of vertices in $V_{1}$
			(resp. $V_{2}$) is joined by a path as described. 
		\end{proof}
		We can now show an effective version of \Cref{prop:Hamiltonian-infinite-3paths}. 
		\begin{proposition}[Computable \Cref{prop:Hamiltonian-infinite-3paths}]
			\label{prop:hamiltonian-3-path-calculable}Let $\G$ be a graph that
			is highly computable, connected, and has either $1$ or $2$ ends.
			Then it admits a bi-infinite Hamiltonian 3-path which is computable.
		\end{proposition}
		
		\begin{proof}
			Let $(v_{i})_{i\in\N}$ be a numbering of the vertex set of the highly
			computable graph $\G$. Now let $f_{0}$ be a 3-path which is bi-extensible
			and visits $v_{0}$. If $\G$ has two ends, then we also require that
			$\G-f_{0}$ has two infinite connected components. In this case we
			do not claim that the path $f_{0}$ can be computed from a description
			of the graph, but it exists and can be specified with finite information.
			After fixing $f_{0}$, we just follow the proof of \Cref{prop:Hamiltonian-infinite-3paths},
			and observe that a sequence of 3-paths $(f_{n})_{n\in\N}$ as in this
			proof can be uniformly computed. That is, there is an algorithm which
			given $n$, computes $f_{n}$. The algorithm proceeds recursively:
			assuming that $(f_{i})_{i\leq n}$ have been computed, we can compute
			$f_{n+1}$ by an exhaustive search. The search is guaranteed to stop,
			and the conditions that we impose on $f_{n+1}$ are decidable thanks
			to Propositions \ref{prop:bi-extensible-decidible-1-end} and \ref{prop:bi-extensible-decidible-2-end}.
			Finally, let $f\colon\Z\to V(\G)$ be the Hamiltonian 3-path on $\G$
			defined by $f(k)=f_{n}(k)$, for $n$ big enough. Then it is clear
			that the computability of $(f_{n})_{n\in\N}$ implies that $f$ is
			computable. 
		\end{proof}
		Now, we are ready to prove \Cref{thm:computable-t1}.
		
		\begin{proof}[Proof of \Cref{thm:computable-t1}]
			Let $\G$ be as in the statement. By \Cref{prop:hamiltonian-3-path-calculable},
			$\G$ admits a bi-infinite Hamiltonian 3-path $f\colon\Z\to V(\G)$
			that is computable. Then it is clear that the translation-like action
			$\ast\colon V(\G)\times\Z\to V(\G)$ defined by $v\ast n=f(f^{-1}(v)+n),\ n\in\Z$,
			is computable.
		\end{proof}
		This readily implies \Cref{thm:computable-translation-like-actions-with-decidable-orbit-problem}
		for groups with one or two ends. 
		\begin{proof}[Proof of \Cref{thm:computable-translation-like-actions-with-decidable-orbit-problem}
			for groups with one or two ends]
			Let $G$ be a finitely generated infinite group with one or two
			ends, and with decidable word problem. Let $S\subset G$ be a finite
			set of generators, and let $\G=\cay(G,S)$ be the associated Cayley
			graph. As $G$ has decidable word problem, this is a highly computable
			graph (\Cref{prop:cayley-graph-of-computable-group-is-computable}).
			Then by \Cref{thm:computable-t1}, $\G$ admits a computable
			and transitive translation-like action by $\Z$. As the vertex set
			of $\G$ is $G$, this is also a computable and transitive translation-like
			action on $G$. This action has decidable orbit membership problem
			for the trivial reason that it has only one orbit. 
		\end{proof}
		\begin{remark}
			As mentioned in the introduction, there is a number of results in
			the theory of infinite graphs that can not have an effective counterpart
			for highly computable graphs. In contrast, we have the following consequences
			of Theorems \ref{thm:t1-translation-like-action-transitivas} and
			\ref{thm:computable-t1}: 
			\begin{enumerate}
				\item A highly computable graph admits a transitive translation-like action
				by $\Z$ if and only if it admits a computable one.
				\item A group with decidable word problem has a Cayley graph with a bi-infinite
				Hamiltonian path if and only if it has a Cayley graph with a computable
				bi-infinite Hamiltonian path. 
				\item The cube of a highly computable graph admits a bi-infinite Hamiltonian
				path if and only if it admits a computable one.
			\end{enumerate}
			The third item should be compared with the following result of D.Bean:
			there is a graph that is highly computable and admits infinite Hamiltonian
			paths, but only uncomputable ones \cite{bean_recursive_1976}. Thus,
			the third item shows that for graphs that are cubes, it is algorithmically
			easier to compute infinite Hamiltonian paths. 
			
			It follows from our results that the problem of \emph{deciding} whether
			a graph admits a bi-infinite Hamiltonian path is also algorithmically
			easier when we restrict ourselves to graphs that are cubes. D.Harel
			proved that the problem of Hamiltonicity is analytic-complete for
			highly computable graphs \cite[Theorem 2]{harel_hamiltonian_1991}.
			On the other hand, it follows from \Cref{thm:t1-translation-like-action-transitivas}
			that for graphs that are cubes, it suffices to check that the graph
			is connected, and has either $1$ or $2$ ends. These conditions are
			undecidable, but are easily seen to be arithmetical \cite{inproceedings}.
			In view of these results, it is natural to ask if these problems are
			easier when we restrict ourselves to graphs that are squares.
		\end{remark}
		
		\begin{question}
			The problem of computing infinite Hamiltonian paths (resp. deciding
			whether an infinite graph is Hamiltonian) on highly computable graphs,
			is easier when we restrict to graphs that are squares?
		\end{question}

		\subsection{Computable normal forms and Stalling's theorem}
		
		In this subsection we prove \Cref{thm:computable-translation-like-actions-with-decidable-orbit-problem}
		for groups with two or more ends. It follows from Stalling's structure
		therem on ends of groups that a group with two or more ends has a
		subgroup isomorphic to $\Z$. We will prove that, if the group has
		solvable word problem, then this subgroup has decidable membership
		problem. This will be obtained from normal forms associated to HNN
		extensions and amalgamated products. We now recall well known facts
		about these constructions, the reader is referred to \cite[Chapter IV]{lyndon_combinatorial_2001}. 
		
		HNN extensions are defined from a group $H=\langle S_{H}|R_{H}\rangle$,
		a symbol $t$ not in $S_{H}$, and an isomorphism $\phi\colon A\to B$
		between subgroups of $H$. The \textit{HNN extension} relative to
		$H$ and $\phi$ is the group with presentation $H\ast_{\phi}=\langle S_{H},t\mid\ R_{H},\ tat^{-1}=\phi(a),\ \forall a\in A\rangle.$
		Now let $T_{A}\subset H$ and $T_{B}\subset H$ be sets of representatives
		for equivalence classes of $H$ modulo $A$ and $B$, respectively.
		The group $H\ast_{\phi}$ admits a normal form associated to the sets
		$T_{A}$ and $T_{B}$. The sequence of group elements $h_{0},t^{\epsilon_{1}},h_{1},\dots,t^{\epsilon_{n}},h_{n}$,
		$\epsilon_{i}\in\{1,-1\}$, is in \textit{normal form} if (1) $h_{0}\in H$,
		(2) if $\epsilon_{i}=-1$, then $h_{i}\in T_{A}$, (3) if $\epsilon_{i}=1$,
		then $h_{i}\in T_{B}$, and (4) there is no subsequence of the form
		$t^{\epsilon},1_{H},t^{-\epsilon}$. For every $g\in H\ast_{\phi}$
		there exists a unique sequence in normal form whose product equals
		$g$ in $H\ast_{\phi}$. 
		
		Amalgamated products are defined from two groups $H=\langle S_{H}|R_{H}\rangle$
		and $K=\langle S_{K}|R_{K}\rangle$, and a group isomorphism $\phi\colon A\to B$,
		with $A\leqslant H$ and $B\leqslant K$. The \textit{amalgamated product} of
		$H$ and $K$ relative to $\phi$, is the group with presentation
		$H\ast_{\phi}K=\langle S_{H},S_{K}|R_{H},R_{K},\ a=\phi(a),\ \forall a\in A\rangle.$
		Now let $T_{A}\subset H$ be a set of representatives for $H$ modulo
		$A$, and let $T_{B}\subset K$ be a set of representatives for $K$
		modulo $B$. The group $H\ast_{\phi}K$ admits a normal form associated
		to the sets $T_{A}$ and $T_{B}$. A sequence of group elements $c_{0},c_{1},\dots,c_{n}$
		is in \textit{normal form }if (1) $c_{0}$ lies in $A$ or $B$, (2)
		$c_{i}$ is in $T_{A}$ or $T_{B}$ for $i\geq1$, (3) $c_{i}\neq1$
		for $i\geq1$, and (4) successive $c_{i}$ alternate between $T_{A}$
		and $T_{B}$. For each element $g\in H\ast_{\phi}K$, there exist
		a unique sequence in normal form whose product equals $g$ in $H\ast_{\phi}K$. 
		
		Stalling's structure theorem relates ends of groups with HNN extensions
		and amalgamated products \cite{dunwoody_cutting_1982}. This result
		asserts that every finitely generated group $G$ with two or more
		ends is either isomorphic to an HNN extension $H\ast_{\phi}$, or
		isomorphic to an amalgamated product $H\ast_{\phi}K$. In both cases,
		the corresponding isomorphism $\phi$ is between finite and proper
		subgroups, and the groups $H$, or $H$ and $K$, are finitely generated
		(see \cite[pages 34 and 43]{cohen_combinatorial_1989}). We will now
		prove that when $G$ has decidable word problem, then the associated
		normal forms are computable. This means that there is an algorithm
		which given a word representing a group element $g$, computes a sequence
		of words such that the corresponding sequence of group elements, is
		a normal form for $g$. The proof is direct, but we were unable to
		find this statement in the literature.
		\begin{proposition}
			\label{prop:stallings-descomposicion-calculable} Let $G$ be a finitely
			generated group with two or more ends and decidable word problem.
			Then the normal form associated to the decomposition of $G$ as HNN
			extension or amalgamated product is computable.
		\end{proposition}
		
		\begin{proof}
			Let us assume first that we are in the first case, so there is a finitely
			generated group $H=\langle S_{H}|R_{H}\rangle$, and an isomorphism
			$\phi\colon A\to B$ between finite subgroups of $H$, such that $G$
			is isomorphic to the HNN extension $H\ast_{\phi}=\langle S_{H},t\mid\ R_{H},\ tat^{-1}=\phi(a),\ a\in A\rangle$.
			A preliminary observation is that $H$ has decidable word problem.
			Indeed, this property is inherited by finitely generated subgroups,
			and $G$ has decidable word problem by hypothesis. The computability
			of the normal form will follow from two simple facts:
			
			First, observe that the finite group $A=\{a_{1},\dots,a_{n}\}$ has
			decidable membership problem in $H$. Indeed, given a word $w\in(S_{H}\cup S_{H}^{-1})^{\ast}$,
			we can decide if $w_{H}\in A$ by checking if $w=_{H}a_{i}$ for $i=1,\dots,m$.
			This is an effective procedure as the word problem of $H$ is decidable,
			and is guaranteed to stop as $A$ is a finite set. As a consequence
			of this, we can also decide if $u_{H}\in Av_{H}$ for any pair of
			words $u,v\in(S_{H}\cup S_{H}^{-1})^{*}$, as this is equivalent to
			decide if $(uv^{-1})_{H}$ lies in $A$. The same is true for $B$.
			
			Second, there is a computably enumerable set $W_{A}\subset(S_{H}\cup S_{H}^{-1})^{*}$
			such that the corresponding set $T_{A}$ of group elements in $H$
			constitute a collection of representatives for $H$ modulo $A$. We
			sketch an algorithm that computably enumerates $W_{A}$ as a computable
			sequence of words. Set $u_{0}$ to be the empty word. Now assume that
			words $u_{0},\dots,u_{n}$ have been selected, and search for a word
			$u_{n+1}\in(S_{H}\cup S_{H}^{-1})^{*}$ such that $(u_{n+1})_{H}$
			does not lie in $A(u_{0})_{H},\dots,A(u_{n})_{H}$. The condition
			that we impose to $u_{n+1}$ is decidable by the observation in the
			previous paragraph, and thus an exhaustive search is guaranteed to
			find a word as required. It is clear that the set $W_{A}$ that we
			obtain is computably enumerable, and that the set $T_{A}$ of the
			group elements of $H$ corresponding to these words is a set of representatives
			for $H$ modulo $A$. A set $W_{B}$ corresponding to $T_{B}$ can
			be enumerated analogously. 
			
			Finally, we note that we can computably enumerate sequences of words
			$w_{0},\dots,w_{n}$ that represent normal forms (with respect to
			$T_{A}$ and $T_{B})$ for all group elements. Indeed, using the fact
			$W_{A}$ and $W_{B}$ are computably enumerable sets, we just have
			to enumerate sequences of words $w_{0},\dots,w_{n}$ such that $w_{0}$
			as an arbitrary element of $(S_{H}\cup S_{H}^{-1})^{*}$, and the
			rest are words from $W_{A},$ $W_{B}$, or $\{t,t^{-1}\}$ that alternate
			as in the definition of normal form. In order to compute the normal
			form of a group element $w_{G}$ given by a word $w$, we just enumerate
			these sequences $w_{1},\dots,w_{n}$ until we find one satisfying
			$w=_{G}w_{1}\dots w_{n}$, this is a decidable question as $G$ has
			decidable word problem. We have proved the computability of normal
			forms as in the statement, in the case where $G$ is a (isomorphic
			to) HNN extension. 
			
			If $G$ is not isomorphic to an HNN extension, then it must be isomorphic
			to an amalgamated product. Then there are two finitely generated groups
			$H=\langle S_{H}|\ R_{H}\rangle$ and $K=\langle S_{K}|R_{K}\rangle$,
			and a group isomorphism $\phi\colon A\to B$, with $A\leqslant H$ and $B\leqslant K$
			finite groups, such that $G$ is isomorphic to $\langle S_{H},S_{K}|R_{H},R_{K},\ a=\phi(a),\ \forall a\in A\rangle$.
			Now the argument is the same as the one given for HNN extensions.
			That is, $A$ and $B$ have decidable membership problem because they
			are finite, and there are two computably enumerable sets of words
			$W_{A}$ and $W_{B}$ corresponding to sets $T_{A}$ and $T_{B}$
			as in the definition of normal form for amalgamated products. This,
			plus the decidability of the word problem, is sufficient to compute
			the normal form of a group element given as a word, by an exhaustive
			search. 
		\end{proof}
		We obtain the following result from the computability of these normal
		forms.
		
		\begin{proposition}
			\label{prop:containment-Z-decidable-membership-problem}Let $G$ be
			a finitely generated group with two or more ends and decidable word
			problem. Then it has a subgroup isomorphic to $\Z$ with decidable
			subgroup membership problem.
		\end{proposition}
		\begin{proof}
			By Stalling's structure theorem, either $G$ is isomorphic to an HNN
			extensions, or $G$ is isomorphic to an amalgamated product. We first
			suppose that $G$ is isomorphic to an HNN extension $H\ast_{\phi}=\langle S_{H},t|\ R_{H},tat^{-1}=\phi(a),a\in A\rangle$.
			Without loss of generality, we will assume that $G$ is equal to this
			group instead of isomorphic, as the decidability of the membership
			problem of an infinite cyclic subgroup is preserved by group isomorphisms.
			We claim that the subgroup of $G$ generated by $t$ has decidable
			membership problem. Indeed, a group element $g$ lies in this subgroup
			if and only if the normal form of $g$ or $g^{-1}$ is $1,t,1\dots,t,1$.
			By \Cref{prop:stallings-descomposicion-calculable}, this normal
			form is computable, and thus we obtain a procedure to decide membership
			in the subgroup of $G$ generated by $t$.
			
			We now consider the case where $G$ is isomorphic to an amalgamated
			product. Then there are two finitely generated groups $H=\langle S_{H}|\ R_{H}\rangle$
			and $K=\langle S_{K}|\ R_{K}\rangle$, and a group isomorphism $\phi:A\to B$,
			with $A\leqslant H$ and $B\leqslant K$ finite groups, such that $G$ is isomorphic
			to $\langle S_{H},S_{K}|R_{H},R_{K},\ a=\phi(a),\ \forall a\in A\rangle$.
			As before, we will assume without loss of generality that $G$ is
			indeed equal to this group. Now let $T_{A}$ and $T_{B}$ be the sets
			defined in \Cref{prop:stallings-descomposicion-calculable} that
			are associated to the computable normal form, and let $u\in T_{A}$,
			$v\in T_{B}$ be both non trivial elements. We claim that the subgroup
			subgroup of $G$ generated by $uv$ is isomorphic to $\Z$, and has
			decidable membership problem. Indeed, a group element $g$ lies in
			this subgroup if and only if the normal form of $g$ or $g^{-1}$
			is $u,v,\dots,u,v$. By \Cref{prop:stallings-descomposicion-calculable},
			this normal form is computable, and thus we obtain a procedure to
			decide membership in the subgroup of $G$ generated by $t$.
		\end{proof}
		We now verify the fact that for translation-like actions coming from
		subgroups, the properties of decidable orbit membership problem and
		decidable subgroup membership problem are equivalent. 
		
		\begin{proposition}
			\label{prop:obvio} Let $H\leqslant G$ be finitely generated groups.
			Then $H$ has decidable membership problem in $G$ if and only if
			the action of $H$ on $G$ by right translations has decidable orbit
			membership problem. 
		\end{proposition}
		
		\begin{proof}
			Let $\ast$ be the action defined by $G\times H\to G$, $(g,h)\mapsto gh$.
			The claim follows from the fact that two elements $g_{1},g_{2}\in G$
			lie in the same $\ast$ orbit if and only if $g_{1}g_{2}{}^{-1}\in H$,
			and an element $g\in G$ lies in $H$ if and only if it lies in the
			same $\ast$ orbit as $1_{G}$.
			
			It is clear how to rewrite this in terms of words, but we fill the
			details for completeness. For the forward implication, let $u,v\in(S\cup S^{-1})^{*}$
			be two words, for which we want to decide whether $u_{G}$, $v_{G}$
			lie in the same orbit. We start by computing the formal inverse of
			$v$, denoted $v^{-1}$, and then check whether the word $uv^{-1}$
			lies in $\{w\in(S\cup S^{-1})^{*}|\ w_{G}\in H\}$. This set is decidable
			for hypothesis. For the reverse implication, assume that the action
			has decidable orbit membership problem. The set $\{w\in(S\cup S^{-1})^{*}|\ w_{G}\in H\}$
			equals the set of words $w\in(S\cup S^{-1})^{*}$ such that $w_{G}$
			and $1_{G}$ lie in the same orbit, which is a decidable set by hypothesis.
			It follows that $H$ has decidable subgroup membership problem in
			$G$. 
		\end{proof}
		We can now finish the proof of \Cref{thm:computable-translation-like-actions-with-decidable-orbit-problem}. 
		\begin{proof}[Proof of \Cref{thm:computable-translation-like-actions-with-decidable-orbit-problem}
			for groups with two or more ends]
			Let $G$ be a finitely generated infinite group with decidable word
			problem and at least two ends. By \Cref{prop:containment-Z-decidable-membership-problem}
			there is an element $c\in G$ such that $\langle c\rangle$ is isomorphic
			to $\Z$, and has decidable subgroup membership problem in $G$. The
			right action $\Z\curvearrowright G$ defined by $g\ast n=gc^{n}$
			has decidable orbit membership problem by \Cref{prop:obvio}. 
			
			It only remains to verify that the function $G\times\Z\to G$, $(g,n)\mapsto g\ast n$
			is computable in the sense of \Cref{subsec:numberings}. This
			is clear, but we write the details for completeness. The group operation
			$f_{1}\colon G\times G\to G$ is computable by \Cref{prop:computabilitiy-groups}.
			Moreover, it is clear that the function $f_{2}\colon\Z\to G$, $n\mapsto c^{n}$
			is computable. Then it follows that the function $f_{3}\colon G\times\Z\to G$,
			$(g,n)\mapsto f_{1}(g,f_{2}(n))$ is computable, being the composition
			of computable functions. But $f_{3}(g,n)=g\ast n$, and thus $\ast$
			is a computable group action. 
		\end{proof}

		\section{Medvedev degrees of effective subshifts}\label{sec:medvedev}
		
		The goal of this section is to prove \Cref{thm:medvedev-degrees-of-effective-subshifts}.
		That is, that on every infinite group with decidable word problem,
		the class of possible Medvedev degrees of effective subshifts is that
		of $\Pi_{1}^{0}$ degrees. 
		
		In summary, our proof is the application of a known construction that
		given a subshift $X\subset A^{\Z}$, outputs a subshift $Y\subset A^{G}$
		whose configurations describe simultaneously translation-like actions
		$\Z\curvearrowright G$, and configurations in $X$. When we require
		this construction to preserve the Medvedev degree of the initial subshift
		$X$, then the existence of a computable translation-like action $\Z\curvearrowright G$
		with decidable orbit membership problem arises as a natural condition. 
		
		Medvedev degrees of subshifts have only been discussed in the literature
		for $G=\Z^{d}$, and for this reason we will review computability
		aspects of the space $A^{G}$ in detail. Given a group $G$ with decidable
		word problem, we will translate computability notions from $A^{\N}$
		to $A^{G}$ using a computable numbering $\nu\colon\N\to G$. We verify
		that the computability notions in this space are independent of the
		chosen numbering, preserved by group isomorphisms, compatible with
		previous notions in the literature \cite{aubrun_notion_2017}, and
		that an effective subshift is the same as a subshift that is effectively
		closed as a set. This equivalence is lost for groups whose word problem
		is algorithmically complex, see \cite{aubrun_notion_2017,barbieri_effective_2024}.
		
		\subsection{Computability notions on the Cantor space }\label{subsec:Computability-on-A^n}
		
		Here we will review some standard concepts from the theory of computability
		on the Cantor space. A modern reference of computability theory on
		uncountable spaces is \cite{brattka_handbook_2021}. 
		
		Let $A$ be a finite alphabet. The set $A^{\N}$ is endowed with the
		pro-discrete topology, for which a sub-basis is the set of cylinders.
		A \textit{cylinder }is a set of the form $[p]=\{x\in A^{\N}\mid x|_{K}=p\}$,
		where $p$ is a \textit{pattern}: a function with from a finite set
		$K\subset\N$ to $A$. We identify a word $w=w_{0}\dots w_{n}\in A^{*}$
		with the pattern $\{0,\dots,n\}\to A$, and thus $[w]=\{x\in A^{\N}|\ x_{0}\dots x_{n}=w_{0}\dots w_{n}\}$. 
		\begin{definition}
			A set $X\subset A^{\N}$ is \textit{effectively closed}, denoted $\Pi_{1}^{0}$,
			if some of the following equivalent conditions hold:
			\begin{enumerate}
				\item The complement of $X$ can be written as $\bigcup_{w\in L}[w]$, for
				a computably enumerable set of words $L\subset A^{*}$.
				\item It is semi-decidable whether a word $w$ satisfies $[w]\cap X=\emptyset$.
				\item It is semi-decidable whether a pattern $p$ satisfies $[p]\cap X=\emptyset$.
			\end{enumerate}
		\end{definition}
		
		\begin{definition}
			A partial function $F\colon D\subset A^{\N}\to B^{\N}$ is \textit{computable}
			when there is a partial computable function on words $f\colon A^{\ast}\to B^{\ast}$
			satisfying the following three conditions:
			\begin{enumerate}
				\item $f$ is monotone for the prefix order on words.
				\item For each $x$ in the domain $D$, the length of $f(x|_{\{0,\dots,k\}})$
				tends to infinity with $k$. 
				\item For every $x$ in the domain $D$, and for every $k\in\N$, there
				is $n$ big enough such that $F(x)(n)$ is the $n$-th letter in the
				word $f(x|_{\{0,\dots,k\}})$. 
			\end{enumerate}
		\end{definition}
		
		It follows from the definition that a computable function must be
		continuous. 
		\begin{example}
			The shift function $\sigma\colon A^{\N}\to A^{\N},\ (\sigma x)(n)=x(n+1)$
			is computable. This is shown by the computable function $s\colon A^{*}\to A^{*}$,
			$s(w_{0}w_{1}\dots w_{n})=w_{1}\dots w_{n}$.
		\end{example}
		
		\begin{definition}
			Let $X\subset A^{\N}$ and $Y\subset B^{\N}$. The sets $X$ and $Y$
			are \textit{computably homeomorphic }if there is a homeomorphism $\Phi\colon X\to Y$
			such that both $\Phi$ and and its inverse are computable functions.
		\end{definition}
		
		\begin{example}
			\label{exa:biycalculable}Let $f\colon\N\to\N$ be a computable bijection,
			and let $F\colon A^{\N}\to A^{\N}$ be defined by $x\mapsto x\circ f$.
			Then $F$ is a computable homeomorphism, and with computable inverse
			$x\mapsto x\circ f^{-1}$.
		\end{example}
		
		\begin{example}
			If $A$ and $B$ are finite alphabets with cardinality at least $2$,
			then the sets $A^{\N}$ and $B^{\N}$ are computably homeomorphic.
			Indeed, the usual homeomorphism between these sets is a computable
			function (see \cite[Theorem 2-97]{hocking_topology_1961}). A simple
			case is when $A=\{0,1,2,3\}$ and $B=\{0,1\}$. Then a computable
			homeomorphism is given by the letter-to-word substitutions $0\mapsto00,\ 1\mapsto01,\ 2\mapsto10,\ 3\mapsto11.$ 
		\end{example}

		\subsection{Medvedev degrees}
		
		Here we review the lattice $\M$ of Medvedev degrees. A survey on
		this topic is \cite{hinman_survey_2012}. 
		\begin{definition}
			Let $X\subset A^{\N}$ and $Y\subset B^{\N}$. We say that $Y$ is
			\textit{Medvedev reducible} to $X$, written $Y\leq_{\M}X$, if there
			is a partial computable function $\Phi$ defined on all elements of
			$X$, and such that $\Phi(X)\subset Y$. We write $X\equiv_{\M}Y$
			when we have both reductions. A \textit{Medvedev degrees} is an equivalence
			class of $\equiv_{\M}$, and we denote by $\M$ the set of Medvedev
			degrees. The pre-order $\leq_{\M}$ becomes a partial order on $\M$,
			and the degree of a set $X$ is denoted by $\deg_{\M}(X)$. 
		\end{definition}
		
		The partially ordered set $(\M,\leq_{\M})$ is indeed a distributive
		lattice with a bottom element $0_{\M}$, and a top element $1_{\M}$.
		We remark that the Medvedev degree of a set $X$ is meaningful when
		we regard $X$ as the set of all solutions to a problem: it measures
		how hard is it to find a solution, where hard means hard to compute.
		For instance, $\deg_{\M}(X)=0_{\M}$ if and only if $X$ has a computable
		point, while $\deg_{\M}(X)=1_{\M}$ if and only if $X$ is empty.
		A prominent sub-lattice of $\M$ is that of $\Pi_{1}^{0}$ degrees:
		\begin{definition}
			A Medvedev degree is called $\Pi_{1}^{0}$ when it is the degree of
			a $\Pi_{1}^{0}$ nonempty subset of $\{0,1\}^{\N}$. 
		\end{definition}

		\subsection{Subshifts}
		
		Here we review standard terminology for subshifts. The reader is referred
		to the book \cite{ceccherini-silberstein_cellular_2018}. 
		
		Let $G$ be a finitely generated group, and let $A$ be a finite alphabet.
		We endow $A^{G}$ with the prodiscrete topology. A \textit{subshift
		}is a subset $X\subset A^{G}$ which is closed and invariant under
		the group action $G\curvearrowright A^{G}$ by left translations $(gx)(h)\mapsto x(g^{-1}h).$
		A \textit{pattern }is a function $p$ from a finite set $K\subset G$
		to $A$, and it determines the\textit{ cylinder $[p]=\{x\in A^{G}\ |\ x|_{K}=p\}.$
		}If $gx\in[p]$ for some $g\in G$, we say that $p$ \textit{appears}
		on $x$. A set of \textit{forbidden }patterns $\F$ defines the subshift
		$X_{\mathcal{\F}}$ of all elements $x\in A^{G}$ where no pattern
		of $\F$ appears in $x$. Every subshift is determined by a maximal
		set of forbidden patterns, but it can have more than one defining
		set of forbidden patterns. A subshift is\textit{ of finite type }(SFT)
		if it can be defined with a finite set of forbidden patterns. 
		
		\subsection{Computability on $A^{G}$  }\label{subsec:Computability-on-A^g}
		
		In this subsection we translate computability notions from $A^{\N}$
		to $A^{G}$, where $G$ is a finitely generated group with decidable
		word problem. Our goal is to provide a definition of Medvedev degree
		for subshifts. In simple words, we will take a computable bijection
		$\nu\colon\N\to G$, and use it to define a homeomorphism from $A^{\N}$
		to $A^{G}$. We declare this homeomorphism to be computable, and in
		this manner we translate to $A^{G}$ the concepts defined in $A^{\N}$.
		This process is well established in the theory of computability on
		uncountable spaces, and is the subject of representation theory. A
		representation plays the same role as a numbering (\Cref{subsec:numberings}),
		but for an uncountable set. 
		
		We recall now some definitions from \cite[Chapter 9]{brattka_handbook_2021}.
		A \textit{represented space} is a pair $(X,\delta)$ where $X$ is
		a set and $\delta$ is a \textit{representation} of $X$: a partial
		surjection $\delta\colon\text{dom}(\delta)\subset A^{\N}\to X$. In
		a represented space $(X,\delta)$, a subset $Y\subset X$ is \textit{effectively
			closed} when $\delta^{-1}(Y)\subset A^{\N}$ is an effectively closed
		set. Moreover, if $(X',\delta'\colon A'{}^{\N}\to X')$ is another
		represented space, a function $F\colon X\to X'$ is \textit{computable}
		when $\delta'{}^{-1}\circ F\circ\delta\colon A^{\N}\to A'{}^{\N}$
		is a computable function. Finally, two representations of the same
		space $X$, $\delta\colon A^{\N}\to X$ and $\delta'\colon A'{}^{\N}\to X$,
		are \textit{equivalent }if the identity function from $(X,\delta)$
		to $(X,\delta')$ is computable. Note that in this case, both representations
		induce the same computability notions on $X$.
		
		In what follows we will focus on a specific representation of $A^{G}$,
		which is also a total function and a homeomorphism. 
		\begin{definition}
			\label{def:gdelta} Let $G$ be a finitely generated group with decidable
			word problem, and let $\nu$ a computable numbering of $G$. We define
			the representation $\delta$ by
			
			\begin{align*}
				\delta\colon A^{\N} & \to A^{G}\\
				x & \mapsto x\circ\nu^{-1}.
			\end{align*}
		\end{definition}
		
		It follows from \Cref{prop:computabilitiy-groups} that a group
		as in the statement admits a computable numbering, and that all these
		numberings are equivalent. In terms of representations, this is expressed
		as follows:
		\begin{proposition}
			\label{prop:equivalent-representations}In \Cref{def:gdelta},
			any two computable numberings induce equivalent representations. 
		\end{proposition}
		
		\begin{proof}
			Let $\nu'$ be another computable numbering of $G$, and let $\delta'$
			be the associated representation of $A^{G}$. Let $F\colon A^{G}\to A^{G}$
			be the identity function. Then $\delta'{}^{-1}\circ F\circ\delta\colon A^{\N}\to A^{\N}$
			is given by $x\mapsto x\circ\nu^{-1}\circ\nu'$. We verify that this
			function is a computable homeomorphism. Indeed, as the numberings
			$\nu$, $\nu'$ are equivalent (\Cref{prop:computabilitiy-groups}),
			the function $\nu^{-1}\circ\nu'\colon\N\to\N$ is a computable bijection
			of $\N$, and this implies that $x\mapsto x\circ\nu^{-1}\circ\nu'$
			is a computable homeomorphism (see \Cref{exa:biycalculable}).
		\end{proof}
		We note that computability notions on $A^{G}$ are also preserved
		by group isomorphisms. 
		\begin{proposition}
			\label{prop:stability}Let $G$ and $G'$ be finitely generated groups
			with decidable word problem, and let $A^{G}$, $A^{G'}$ be endowed
			with the representation in \Cref{def:gdelta}. If $f\colon G\to G'$
			is a group isomorphism, then the associated function $F\colon A^{G'}\to A^{G}$,
			$x\mapsto x\circ f$ is a computable homeomorphism.
		\end{proposition}
		
		The proof is similar to the proof of \Cref{prop:equivalent-representations},
		but applying the third item in \Cref{prop:computabilitiy-groups}.
		This means that the computability notions on $A^{G}$ are preserved
		if we rename group elements (for example, by taking different presentations
		of the same group). We are ready to define the Medvedev degree of
		a subset of $A^{G}$. 
		\begin{definition}
			Let $G$ be a finitely generated group with decidable word problem.
			Given a subset $X\subset A^{G}$, we define $\deg_{\M}(X)=\deg_{\M}(\delta^{-1}X)$.
		\end{definition}
		
		This definition does not depend on $\delta$, as long as $\delta$
		comes from a computable numbering of $G$. We now turn our attention
		to effectively closed subsets of $A^{G}$, and subshifts. 
		\begin{proposition}
			\label{prop:ef-cerrado}Let $G$ be a finitely generated group with
			decidable word problem. Then a subset $X\subset A^{G}$ is effectively
			closed if and only if it is semi-decidable whether a pattern $p\colon K\subset G\to A$
			satisfies $[p]\cap X=\emptyset$.
		\end{proposition}
		
		\begin{proof}
			We only prove the forward implication, the converse being similar.
			Given a pattern $p\colon K\subset G\to A$, we start by computing
			a pattern $p'\colon K'\subset\N\to A$ such that $p=p'\circ\nu$.
			Then $[p]\cap X=\emptyset$ if and only if $[p']\cap\delta^{-1}(X)=\emptyset$.
			But the latter relation is semi-decidable on $p'$ as $\delta^{-1}(X)$
			is effectively closed in $A^{\N}$. 
		\end{proof}
		In \cite{aubrun_notion_2017} the authors introduced a notion of effectiveness
		for subshift on general finitely generated groups. This notion is
		not explicitely associated to $A^{G}$ as a represented space (or
		a computable metric space), but we shall verify now that for groups
		with decidable word problem, these approaches are equivalent. 
		
		The following definitions are taken from \cite{aubrun_notion_2017}.
		A \textit{pattern coding} $c$ is a finite set of tuples $\{(w_{1},a_{1}),\dots,(w_{k},a_{k})\}$,
		where $w_{i}\in S^{*}$ and $a_{i}\in A$, and is \textit{consistent}
		when $w_{i}=_{G}w_{j}$ implies $a_{i}=a_{j}$. A consistent pattern
		coding defines a pattern $p(c)\colon K\subset G\to A$, where $K$
		equals $\{(w_{1})_{G},\dots,(w_{k})_{G}\}$, and $p((w_{i})_{G})=a_{i}$.
		A set of pattern codings $\mathcal{C}$ defines the subshift $X_{\mathcal{C}}$
		of all elements $x\in A^{G}$ such that no pattern of the form $p(c)$
		appears in $x$, where $c$ ranges over $\mathcal{C}$. A subshift
		$X$ is \textit{effective }if there is a computably enumerable set
		of pattern codings $\mathcal{C}$ such that $X=X_{\mathcal{C}}$. 
		\begin{proposition}
			\label{prop:characterization-effective-subshifts-1}Let $G$ be a
			finitely generated group with decidable word problem. Then a subshift
			$X\subset A^{G}$ is effective if and only if it is an effectively
			closed subset of $A^{G}$.
		\end{proposition}
		
		\begin{proof}
			If a subshift is an effectively closed subset of $A^{G}$, then by
			\Cref{prop:ef-cerrado} the set of all patterns $p$ with $[p]\cap X=\emptyset$
			is computably enumerable. Let $\mathcal{F}$ be this set of patterns,
			and let $\mathcal{C}$ be the set of all pattern codings associated
			to patterns in $\mathcal{F}$. It is clear that $\mathcal{C}$ is
			computably enumerable and $X=X_{\mathcal{C}}$, so $X$ is an effective
			subshift as well.
			
			We now consider the other direction. In \cite[Lemma 2.3]{aubrun_notion_2017}
			it is shown that for a recursively presented group and in particular
			one with decidable word problem, an effective subshift has a maximal
			-for inclusion- computably enumerable set of pattern codings associated
			to forbidden patterns. Given an effective subshift $X$, we can write
			$X=X_{\mathcal{C}}$, where $\mathcal{C}$ is a maximal -for inclusion-
			set of defining forbidden pattern codings. As $G$ has decidable word
			problem, the set of consistent pattern codings is decidable, and thus
			we can computably discard those pattern codings that are not consistent.
			This, plus the previous fact, proves that the set of all patterns
			$p$ with $[p]\cap X=\emptyset$ is computably enumerable. By \Cref{prop:ef-cerrado},
			it follows that the set $X$ is effectively closed.
		\end{proof}

		Let us now make some comments about the computability of the action
		$G\curvearrowright A^{G}$ by translations. It follows from \Cref{prop:computabilitiy-groups}
		and \Cref{exa:biycalculable} that this action is computable. 
		\begin{proposition}
			\label{prop:accion-calculable}Let $G$ be a finitely generated group
			with decidable word problem. Then the group action $G\curvearrowright A^{G}$
			is computable. 
		\end{proposition}
		
		It follows from \Cref{prop:computabilitiy-groups} that all numberings
		of a group $G$ as above that make the left (resp. right) action $G\curvearrowright G$,
		$(g,h)\mapsto gh$ computable, are equivalent. In other word, the
		action of a group on itself characterizes those computable numberings
		of the group. This is a well studied subject, and leads to the notion
		of \emph{computable dimension }of a group. See for instance \cite{goncharov_computable_2003}.
		It is natural then to ask whether something analogous happens for
		representations of the space $A^{G}$:
		\begin{question}
			Let $G$ be a finitely generated group with decidable word problem.
			Are all representations of the space $A^{G}$ that make the action
			$G\curvearrowright A^{G}$ computable equivalent?
		\end{question}

		\subsection{The subshift of translation-like actions by $\Z$, and the
			proof of \Cref{thm:medvedev-degrees-of-effective-subshifts}}
		
		In this subsection we finally prove \Cref{thm:medvedev-degrees-of-effective-subshifts}.
		Our standing assumption is that $G$ is a finitely generated group,
		$S\subset G$ is a finite set of generators, and $J\in\N$. When we
		need $G$ to have decidable word problem, we will specify it. 
		\begin{definition}
			We define $T_{J}(\Z,G)$ as the set of all translation-like actions
			$\ast\colon G\times\Z\to G$, such that $\ensuremath{\{d_{S}(g,g\ast1)\mid g\in G}\}$
			is bounded by $J$.
		\end{definition}
		
		Consider now the finite alphabet $B=B(1_{G},J)\times B(1_{G},J)$,
		where $B(1_{G},J)$ is the ball $\{g\in G\mid d_{S}(g,1_{G})\leq J\}$.
		Every translation-like action $\ast\in T_{J}(\Z,G)$ defines a configuration
		in $B^{G}$, denoted $x_{\ast}$, by the condition 
		\[
		\forall g\in G\hspace{1em}x_{\ast}(g)=(l,r)\iff g\ast-1=gl\text{ and }g\ast1=gr.
		\]
		
		\begin{definition}
			We define $X_{J}(\Z,G)$ as the set $\{x_{\ast}\in B^{G}\mid\ast\in T_{J}(\Z,G)\}$. 
		\end{definition}
		
		\begin{figure}
			\begin{center}\includegraphics[width=0.6\columnwidth]{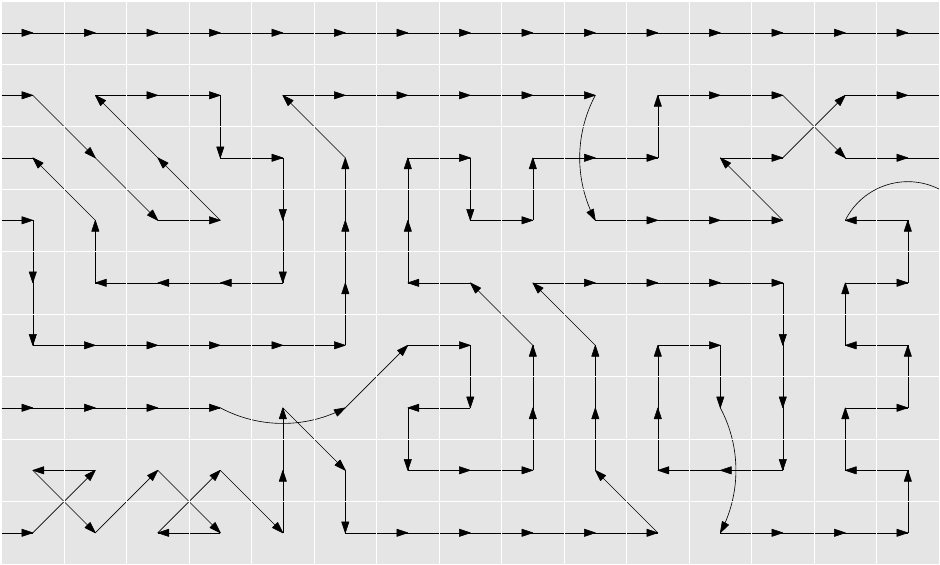}\end{center}
			
			\caption{Representation of some orbits of a translation-like action in $T_{2}(\Z,\Z^{2})$,
				or alternatively, a finite pattern in a configuration in $X_{2}(\Z,\Z^{2})$.
				In this case, $\Z^{2}$ is endowed with the set of four generators
				$S=\{(\pm1,0),(0,\pm1)\}$.}\label{fig:subshift-translation-like-action}
		\end{figure}
		The informal idea is to interpret $x(g)=(l,r)$ as a pair of arrows:
		$g$ has an outgoing arrow to $gr$, and an incoming arrow from $gl$. See \Cref{fig:subshift-translation-like-action}.
		\begin{proposition}
			\label{prop:X_j-es-subshift} The set $X_{J}(\Z,G)$ is a subshift.
			If $G$ has decidable word problem, then it is an effective subshift. 
		\end{proposition}
		
\begin{figure}
			\begin{center}\includegraphics[width=0.6\columnwidth]{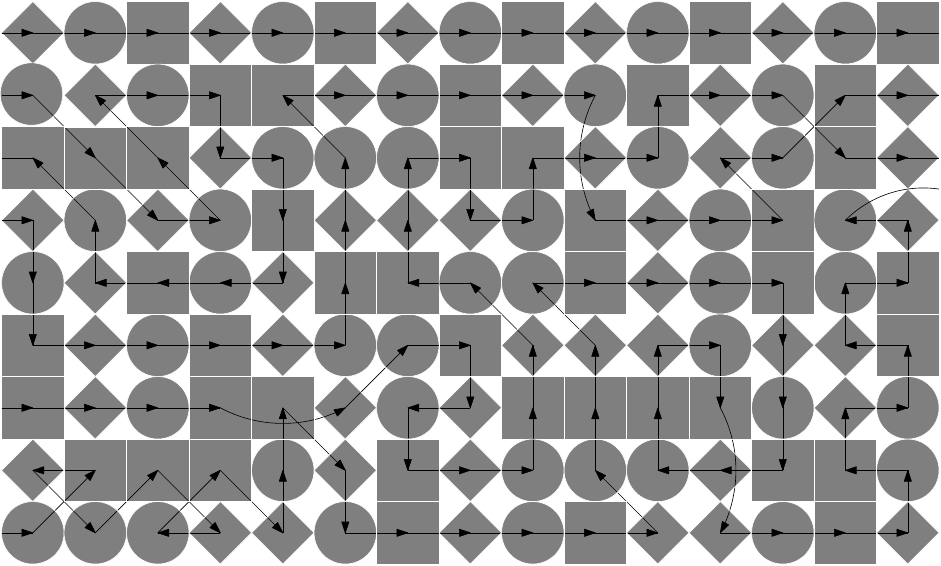}\end{center}
			
			\label{fig:translation-like-actions-y-subshift}\caption{Representation of a finite pattern in a of configuration in $Y[X_{2}(\Z,\Z^{2})]$.
				Here $A$ is the alphabet $\{\text{circle},\text{square},\text{rhombus}\}$,
				and $Y\subset A^{\Z}$ is the subshift of all sequences that
				alternate circle, square, and rhombus in that order.}
		\end{figure}

		\begin{proof}
			We define for each element $x\in B^{G}$ a function $\ast_{x}\colon G\times\Z\to G$,
			which may not be a group action. $L$ and $R$ stand for the projections
			$B\to B(1_{G},J)$ to the left and right coordinate, respectively.
			For $m\in\Z_{\geq0}$ and $g\in G$, define $g\ast_{x}m$ by setting
			$g\ast_{x}0=g$, $g\ast_{x}1=gR(x(g))$, and $g\ast_{x}(m+1)=(g\ast_{x}m)\ast_{x}1.$
			For $m\in\Z_{\leq0}$, define $g\ast_{x}m$ by $g\ast_{x}-1=gL(x(g))$
			and $g\ast_{x}(m-1)=(g\ast_{x}m)\ast_{x}-1$. 
			
			If $p\colon K\subset G\to B$ is a pattern and $m\in\Z$, we give
			to $g\ast_{p}m$ the same meaning as before, as long as it is defined.
			Note that for arbitrary $x\in B^{G}$ and $n,m\in\Z$, the relation
			$(g\ast_{x}n)\ast_{x}m=g\ast_{x}(n+m)$ is not guaranteed to hold,
			but it does hold when $n$ and $m$ have the same sign. 
			
			Let $\J$ be the set of all patterns $p\colon B(1_{G},n)\to B$, $n\in\N$,
			such that some of the following conditions occur:
			\begin{enumerate}
				\item $(1_{G}\ast_{p}1)\ast_{p}-1\neq1_{G}$.
				\item $(1_{G}\ast_{p}-1)\ast_{p}1\neq1_{G}$ .
				\item For some $m\in\Z-\{0\}$, $1_{G}\ast_{p}m=1_{G}$.
			\end{enumerate}
			We claim that $X_{J}(\Z,G)=X_{\J}$. The inclusion $X_{J}(\Z,G)\subset X_{\J}$
			is striaghtforward. Indeed, given $\ast\in T_{J}(\Z,G)$, it is clear
			that no pattern of $\J$ may appear on $x_{\ast}$ by the definition
			of group action and translation-like action. 
			
			We prove now that $X_{\J}\subset X_{J}(\Z,G)$. Let $x\in X_{\J}$
			be an arbitrary element. We first prove that $\ast_{x}$ is a translation-like
			action in $T_{J}(\Z,G)$. Indeed, it follows from the forbidden patterns
			in $\J$ that for every $g\in G$, $(g\ast_{x}1)\ast_{x}-1=(g\ast_{x}-1)\ast_{x}1=g$.
			Then an easy induction on $\max\{|n|,|m|\}$ shows that $(g\ast_{x}n)\ast_{x}m=g\ast_{x}(n+m)$
			for all $n,m\in\Z$. Thus $\ast_{x}$ is a group action. This action
			is free by the third condition on the set $\J$ foribdden patterns,
			and the boundedness condition comes from the alphabet chosen. Thus
			$\ast_{x}$ is a translation-like action in $T_{J}(\Z,G)$, and then
			$x_{(\ast_{x})}$ lies in $X_{J}(\Z,G)$ by definition. But $x=x_{(\ast_{x})}$,
			so it follows that $x$ lies in $X_{J}(\Z,G)$. As $x$ was an arbitrary
			element from $X_{\J}$, we obtain the desired inclusion $X_{\J}\subset X_{J}(\Z,G)$. 
			
			We now verify that, having $G$ decidable word problem, the subshift
			$X_{J}(\Z,G)$ is effective. The definition of $\ast_{p}$ above is
			recursive: given a pattern $p$ on alphabet $B$ and $m\in\Z$, we
			can decide if the group element $1_{G}\ast_{p}m$ is defined, and
			compute it. This shows that the conditions on patterns (1), (2), and
			(3) are decidable over patterns, and thus that $\J$ is a decidable
			set of patterns. Thus $X_{\J}$ is an effective subshift. 
		\end{proof}
		We now describe a subshift on $G$ whose elements describe, simultaneously,
		translation-like actions, and configurations from a subshift over
		$\Z$. See \Cref{fig:translation-like-actions-y-subshift}. Let
		$A$ be an arbitrary finite alphabet, and let $B$ be the alphabet
		already defined and which depends on the natural number $J$. Elements
		of $(A\times B)^{G}$ can be conveniently written as $(y,x)$ for
		$y\in A^{G}$ and $x\in B^{G}$. We will write $\pi_{A}\colon A\times B\to A$
		and $\pi_{B}\colon A\times B\to B$ for the projections to the first
		and second coordinate, respectively.
		\begin{definition}
			For a one dimensional subshift $Y\subset A^{\Z}$, we define $Y[X_{J}(\Z,G)]$
			as the set of all configurations $(y,x)\in(A\times B)^{G}$ such that
			the following two conditions are satisfied:
			\begin{enumerate}
				\item $x\in X_{J}(\Z,G)$.
				\item For every $g\in G$, the element $m\mapsto y(g\ast_{x}m)$, $m\in\Z$,
				lies in $Y$.
			\end{enumerate}
		\end{definition}
		
		\begin{proposition}
			\label{prop:Y=00005BX_j=00005D-es-subshift}The set $Y[X_{J}(\Z,G)]$
			is a subshift. If $G$ has decidable word problem and $Y$ is an effective
			subshift, then $Y[X_{J}(\Z,G)]$ is an effective subshift. 
		\end{proposition}
		
		\begin{proof}
			Let $\F$ be the set of all patterns in $\Z$ that do not appear in
			$X$, so that $X=X_{\F}$, and let $\J$ be as in the proof of \Cref{prop:X_j-es-subshift}.
			That is, $X_{\J}=X_{J}(\Z,G)$. Define $\mathcal{H}$ to be the set
			of all patterns $p\colon B(1_{G},n)\to A\times B$, $n\in\N$, such
			that some of the following conditions hold: 
			\begin{enumerate}
				\item The pattern $\pi_{B}\circ p\colon B(1_{G},n)\to B$ lies in $\J$.
				\item Let $q=\pi_{B}\circ p\colon B(1_{G},n)\to B$. For some $m\in\N$,
				the elements $g\ast_{q}1,\dots,g\ast_{q}m$ are all defined, lie in
				$\text{\ensuremath{B(1_{G},n)}},$ and the pattern $r\colon\{1,\dots,m\}\subset\Z\to A,r(k)=\pi_{A}(g\ast_{q}k)$
				lies in $\F$.
			\end{enumerate}
			It is a rutinary verification that $x\in Y[X_{J}(\Z,G)]$ if and only
			if $x\in X_{\mathcal{H}}$. This shows that $Y[X_{J}(\Z,G)]$ is a
			subshift. 
			
			Now assume that $G$ has decidable word problem, and $\mathcal{F}$
			is a computably enumerable set. Then the first condition of $\mathcal{H}$
			is decidable on patterns: given a pattern $p$, we can compute the
			pattern $q=\pi_{B}\circ p\colon B(1_{G},n)\to B$, and we already
			proved that $\J$ is a decidable set. The second condition of $\mathcal{H}$
			is semi-decidable: given $p$ and $m\in\N$, we can compute the pattern
			$r$, and semi-decide whether it lies in $\F$. It follows that $\mathcal{H}$
			is a computably enumerable set. 
		\end{proof}
		These constructions were introduced in \cite{jeandel_translationlike_2015}
		in the more general case where there is a finitely generated goup
		$H$ instead of $\Z$. We will write $X_{J}(H,G)$ and $Y[X_{J}(\Z,G)]$
		with the same meaning as before, but only for reference purposes.
		It is natural to ask what properties are preserved by the map $Y\mapsto Y[X_{J}(H,G)]$
		that sends a subshift on $H$ to a subshift on $G$. The following
		is known:
		\begin{enumerate}
			\item In \cite{jeandel_translationlike_2015}, E. Jeandel proved that when
			$H$ is a finitely presented group, this map preserves weak aperiodicity,
			and the property of being empty/nonempty. This was used to show the
			existence of weakly aperiodic subshifts on new groups, and the undecidability
			of the emptiness problem for subshifts of finite type on new groups. 
			\item In \cite{barbieri_entropies_2021}, S. Barbieri proved that when $H$
			and $G$ are amenable groups, the topological entropy $h$ satisfies
			the formula $h(Y[X_{J}(H,G)])=h(Y)+h(X_{J}(H,G)).$ This was used
			to classify the entropy of subshifts of finite type on some amenable
			groups.
		\end{enumerate}
		In the present paper we are interested in the algorithmic complexity
		of subshifts. We already verified that $Y\mapsto Y[X_{J}(\Z,G)]$
		preserves the property of being an effective subshift, which is folklore.
		In the following result we use \Cref{thm:computable-translation-like-actions-with-decidable-orbit-problem}
		to show that $Y\mapsto Y[X_{J}(\Z,G)]$ also preserves the Medvedev
		degree of a subshift when $J$ is big enough.
		\begin{theorem}
			\label{thm:preservaM}Let $G$ be a finitely generated infinite group
			with decidable word problem, and suppose that $J$ is big enough so
			that $T_{J}(\Z,G)$ contains an element as in \Cref{thm:computable-translation-like-actions-with-decidable-orbit-problem}.
			Then for every subshift $Y\subset A^{\Z}$,
			\[
			Y\equiv_{\M}Y[X_{J}(\Z,G)].
			\]
		\end{theorem}
		
		\begin{proof}
			Recall that we have a Medvedev reduction $Y\geq_{\M}X$ when there
			a computable function $\Phi$ defined on all elements of $Y$, and
			with $\Phi(Y)\subset X$. Intuitively, this means that there is an
			algorithm which from any element in $Y$, is able to compute an element
			in $X$. In our case, we will consider computable functions between
			the represented spaces $A^{\Z}$ and $(A\times B)^{G}$, in the sense
			of \Cref{subsec:Computability-on-A^g}. 
			
			Let $\ast$ be a translation-like action as in \Cref{thm:computable-translation-like-actions-with-decidable-orbit-problem},
			and let $J$ be big enough so that $\ast$ lies in $T_{J}(\Z,G)$.
			We first prove the inequality $Y\geq_{\M}Y[X_{J}(\Z,G)]$. The intuitive
			idea is as follows. Given an element $y\in Y$, we define an element
			$(z,x)\in Y[X_{J}(\Z,G)]$ by setting $x=x_{\ast}$ (a computabe point
			of $B^{G}$ because $\ast$ is a computable function), and on each
			orbit described by $x_{\ast}$, we copy the sequence $y$. The fact
			that $\ast$ has decidable orbit membership problem is fundamental:
			when we compute the new element $z\in A^{G}$, we need to know if
			two arbitrary group elements $g$, $h$ can be colored independently
			(when they lie in different orbits by $\ast$), or the color of one
			of them determines the color of the other (when they lie in the same
			orbit by $\ast$). 
			
			Let $(g_{n})_{n\in\N}$ be a computable numbering of $G$. We compute
			a set of representatives for orbits of $\ast$ as follows. Define
			a decidable set $I\subset\N$ by the condition that $n\in I$ when
			$g_{n}$ is the first element in its own orbit that appears in the
			numbering $(g_{n})_{n\in\N}$. This condition is decidable because
			$\ast$ has decidable orbit membership problem. Thus $\{g_{i}\mid i\in I\}$
			contains exactly one representative for each orbit of $\ast$. 
			
			We now define a computable function $\Psi_{A}\colon A^{\Z}\to A^{G}$
			as follows. On input $y$, we define $\Psi_{A}(y)$ by the expression
			\[
			\Psi_{A}(y)(g_{i}\ast n)=y(n),\hspace{1em}i\in I,\ n\in\Z.
			\]
			The sets $\{g_{i}\ast n\mid n\in\Z\}$ partition $G$ when we range
			$i\in I$, and thus we defined $\Psi_{A}(y)(g)$ for all $g\in G$.
			To see that $\Psi_{A}$ is a computable function, we exhibit a procedure
			that given $y\in A^{\Z}$ and $g\in G$, outputs $\Psi_{A}(y)(g)$.
			First, compute $i\in I$ such that $g$ lies in the same orbit as
			$g_{i}$. This is possible as $I$ is a decidable set, and $\ast$
			has decidable orbit membership problem. Then we use the fact that
			the action $\ast$ is computable to find $n\in\Z$ satisfying $g=g_{i}\ast n$.
			Finally, output $y(n)$. As mentioned, this proves that $\Psi_{A}$
			is a computable function. 
			
			Let $\Psi_{B}\colon A^{\Z}\to B^{G}$ be the function with constant
			value $x_{\ast}$, which is a computable because $x_{\ast}$ is a
			computable point. We define now a function $\Psi\colon A^{\Z}\to(A\times B)^{G}$
			by $z\mapsto\Psi(z)=(\Psi_{A}(z),\Psi_{B}(z))$. The function $\Psi$
			is clearly computable, and we have $\Psi(Y)\subset Y[X_{J}(\Z,G)]$
			by construction. This proves the desired inequality $Y\geq_{\M}Y[X_{J}(\Z,G)]$.
			
			The remaining inequality $Y[X_{J}(\Z,G)]\geq_{\M}Y$ is clear. From
			any element in $Y[X_{J}(\Z,G)]$ we can compute an element in $Y$:
			on input $(z,x)$ we just have to follow the arrows from $1_{G}$,
			read the $A$ component of the alphabet, and the sequence obtained
			lies in $Y$. More formally, we define the function $\Phi\colon Y[X_{J}(\Z,G)]\to Y$
			by the expression 
			\[
			\Phi(z,x)(n)=z(1_{G}\ast n),\ n\in\Z.
			\]
			It is clear from the expression above that $\Phi$ is a computable
			function. This proves the desired inequality $Y\geq_{\M}Y[X_{J}(\Z,G)]$.
		\end{proof}
		We are now ready to prove \Cref{thm:medvedev-degrees-of-effective-subshifts}.
		\begin{proof}[Proof of \Cref{thm:medvedev-degrees-of-effective-subshifts}]
			By \Cref{prop:characterization-effective-subshifts-1}, the
			Medvedev degree of every effective subshift on $G$ is a $\Pi_{1}^{0}$
			degree. It follows that the class of Medvedev degrees of effective
			subshifts on $G$ is contained in the class of $\Pi_{1}^{0}$ Medvedev
			degrees.
			
			We now prove that every $\Pi_{1}^{0}$ Medvedev degree is attained
			by a subshift on $G$. Let $P\subset\{0,1\}^{\N}$ be an effectively
			closed set. By Miller's theorem \cite[Proposition 3.1]{miller_two_2012},
			there is an effective subshift on $Y$ on $\Z$, such that $P\equiv_{\M}Y$.
			Suppose that $J$ is big enough so that $T_{J}(\Z,G)$ contains an
			element as in \Cref{thm:computable-translation-like-actions-with-decidable-orbit-problem}.
			Then the subshift $Y[X_{J}(\Z,G)]$ is effective by \Cref{prop:X_j-es-subshift},
			and it satisfies $P\equiv_{\M}Y[X_{J}(\Z,G)]$ by \Cref{thm:preservaM}.
			This finishes the proof.
		\end{proof}
		Our proof of \Cref{thm:medvedev-degrees-of-effective-subshifts}
		has made extensive use of the hypothesis of decidable word problem,
		and it is unclear whether a similar method could work for recursively
		presented groups. 
		\begin{question}
			Let $G$ be a recursively presented infinite group. Is it true that
			effective subshifts on $G$ attain all $\Pi_{1}^{0}$ Medvedev degrees?
		\end{question}
		
		Despite we have not considered recursively presented groups here,
		it can be proved that for recursively presented groups, the Medvedev
		degree of an effective subshift must be a $\Pi_{1}^{0}$ degree \cite[Section 3]{barbieri_effective_2024}.

	\bibliographystyle{abbrv}
	\bibliography{main.bbl}
\end{document}